\documentclass[12pt,a4paper,reqno]{amsart}
\usepackage[utf8]{inputenc}
\usepackage{amssymb,amsmath,amsfonts,eucal,mathrsfs,amsthm,enumerate} 
\usepackage[colorinlistoftodos]{todonotes}
\usepackage[allcolors=blue]{hyperref}
\usepackage{empheq}
\usepackage{amsrefs}

\usepackage[top=5.8cm,bottom=5.8cm,left=4.1cm,right=4.1cm]{geometry}

\usepackage{hyperref}
\hypersetup{
    colorlinks=true,
    filecolor=magenta,      
    linkcolor={blue!50!black},
    citecolor={red!50!black},
    urlcolor={blue!80!black}
}

\newtheorem{theorem}{Theorem}
\newtheorem{proposition}[theorem]{Proposition}
\newtheorem{lemma}[theorem]{Lemma}

\newtheorem{corollary}[theorem]{Corollary}

\theoremstyle{definition}

\newcommand{\CU}{\mathcal{U}}
\newcommand{\Vol}{\mathrm{Vol}}
\newcommand{\R}{\mathbb{R}}
\newcommand{\Fi}{\mathbb{F}}
\newcommand{\Q}{\mathbb{Q}}
\newcommand{\Sf}{\mathbb{S}}
\newcommand{\C}{\mathbb{C}}
\newcommand{\Hy}{\mathbb{H}}

\newcommand{\Ric}{\mbox{Ric }}

\newcommand{\End}{\mbox{End}}
\newcommand{\Hom}{\mbox{Hom}}

\newcommand{\ind}{{\rm Index}\, }
\newcommand{\Tp}{T_A^{[p]}}
\newcommand{\tr}{\mathrm{tr}}

\def\span{{\rm{span}}}
\def\bea{\begin{eqnarray*} }
\def\eea{\end{eqnarray*} }
\def\CP{\mathord{\mathbb C}\mathord{\mathbb P}}
\def\HP{\mathord{\mathbb H}\mathord{\mathbb P}}

\def\beq{\begin{equation}}
\def\Z{\mathord{\mathbb Z}}
\def\N{\mathord{\mathbb N}}
\def\B{\mathcal{B}}
\def\<{{\langle}}
\def\>{{\rangle}}

\def\n{\nabla}

\def\a{\alpha}

\def\be{\begin{equation} }
\def\ee{\end{equation} }

\def\proof{\noindent{\it Proof:  }}
\def\qed{\ifhmode\unskip\nobreak\fi\ifmmode\ifinner
\else\hskip5 pt \fi\fi\hbox{\hskip5 pt \vrule width4 pt
height6 pt  depth1.5 pt \hskip 1pt }}

\begin{document}

\title{Homology vanishing theorems for pinched submanifolds}
\author{Christos-Raent Onti}
\address{Department of Mathematics and Statistics, University of Cyprus, 1678 Nicosia, Cyprus}
\email{onti.christos-raent@ucy.ac.cy}
\author{Theodoros Vlachos}
\address{Department of Mathematics, University of Ioannina, 45110 Ioannina, Greece}
\email{tvlachos@uoi.gr}
\date{}
\maketitle

\renewcommand{\thefootnote}{\fnsymbol{footnote}} 
\renewcommand{\thefootnote}{\arabic{footnote}} 

\renewcommand{\thefootnote}{\fnsymbol{footnote}} 
\renewcommand{\thefootnote}{\arabic{footnote}} 

\begin{abstract}
We investigate the geometry and topology of submanifolds under a sharp pinching condition involving 
extrinsic invariants like the mean curvature and the length of the second fundamental form. Homology 
vanishing results are given that strengthen and sharpen previous ones. In addition, an integral bound is 
provided for the Bochner operator of compact Euclidean submanifolds in terms of the Betti numbers.
\end{abstract}

\section{Introduction}

A fundamental problem in differential geometry is to investigate the relationship between geometry and 
topology of Riemannian manifolds. The same question can be raised from  the point of view of submanifold 
geometry. Indeed, it has been an active field of research to study the effect of pinching conditions, involving intrinsic and 
extrinsic curvature invariants, on the geometry or the topology of submanifolds. 
An important result in this direction was given by Simons  \cite{Simons68} for minimal submanifolds of spheres. Since then, plenty of geometric and topological rigidity results  have been obtained, under various pinching conditions, see for instance \cite{adc, an, CDCK70, fuxu08, guxu12,  guxu13, Savo14, ybshen,  gauchman84, Vlachos02, hv, Vlachos04}. 
In particular, Lawson and Simons \cite{LS73} showed that certain 
bounds on the second fundamental form for submanifolds of spheres force homology 
groups to vanish. In their approach, the second variation of area was exploited  to rule out stable minimal currents in
certain dimensions, and since one can minimize area in a homology class, this
trivializes integral homology.

The aim of the present paper is to study the geometry and topology of submanifolds under a sharp pinching condition 
involving the mean curvature and the length of the second fundamental form. Throughout the paper, $S$ denotes the 
{\it squared length} of the second fundamental form $\a_f$ of an isometric 
immersion $f$, while the \textit{mean curvature} is defined as the length $H=\| \mathcal H \|$ of the 
\textit{mean curvature vector field} given by $\mathcal H=(\mathrm{tr}\, \alpha_f)/n$, where $\mathrm{tr}$ 
means taking the trace. The choice of the pinching condition is  
inspired by the standard immersion of a torus 
$$
\mathbb T^n_p(r)=\mathbb{S}^p(r) \times \mathbb{S}^{n-p}(\sqrt{1-r^2})
$$
into the unit sphere $\mathbb{S}^{n+1}$, where $\mathbb{S}^p(r)$ denotes the $p$-dimensional sphere 
of radius $r<1$. The principal curvatures are $\sqrt{1-r^2}/r$ and $-r/\sqrt{1-r^2}$ of multiplicity $p$ and 
$n-p$, respectively. A direct computation gives that  $S=a(n,p,H,1)$ if $r^2\geq p/n$ whereas $S>a(n,p,H,1)$ if $r^2< p/n$, where 
\begin{equation*}
a(n,p,t,c)=nc+\frac{n^3t^2}{2p(n-p)}-\frac{n|n-2p|t}{2p(n-p)}\sqrt{
n^2t^2+4p(n-p)c},\; t,c\in \mathbb{R}.
\end{equation*}

This example served as the motivation to study isometric immersions $f\colon M^n\to \tilde M^m$ between 
Riemannian manifolds that satisfy the pinching condition
\be\label{1}
S\leq a(n,p,H,c) 
\ee
pointwise, where $p$ is an integer with $1\leq p\leq n/2$, under the mild assumption 
that the curvature operator of the ambient manifold $\tilde M^m$ is bounded from below by a constant $c$. 
In particular, this includes the case of submanifolds of spaces of constant sectional curvature.
Submanifolds that satisfy condition \eqref{1} have already been investigated by several authors  under geometric assumptions (see for instance \cite{adc, Savo14, lin15, ybshen, Xu1993, CN, WX, Zhang, ZZ}). In addition, the topology of submanifolds
that satisfy strict inequality in \eqref{1} was studied in \cite{leung83, sx97.2, Vlachos07, xugu10, guxu12} without any geometric assumptions.

In the present paper, we investigate the geometry and topology of submanifolds that satisfy the 
pinching condition \eqref{1}. We show that this 
pinching condition either forces homology to vanish in a range of intermediate dimensions, or completely 
determines the pinched submanifold. Our results are  sharp and strengthen previous ones (see \cite{an, CS19, fuxu08, Gr, lin15, leung83, Savo14, ybshen, sx97.2, Vlachos07}). A key ingredient in our approach is the Bochner technique that is based 
on the Bochner-Weitzenb\"ock formula. This states that the Laplacian of every $p$-form 
$\omega\in\Omega^p(M^n)$ on a manifold $M^n$ is given by
\be\label{boch.form}
\Delta\omega=\n^*\n\omega+\B^{[p]}\omega,
\ee
where $\n^*\n$ is 
the rough Laplacian and 
$$
\B^{[p]}\colon \Omega^p(M^n)\to \Omega^p(M^n)
$$
is a certain symmetric endomorphism of the bundle of $p$-forms, called the {\it Bochner operator}.
The Bochner-Weitzenb\"ock formula implies that every harmonic $p$-form vanishes on compact manifolds, 
provided that $\B^{[p]}$ is positive definite. Since every de Rham cohomology class is represented by a 
harmonic form, it follows that the $p$-th de Rham cohomology group $H^p(M^n; \R)$ vanishes.

The difficulty when dealing with the Bochner operator $\B^{[p]}$ arises from the fact that it is unmanageable 
for $p>1$. Our approach is inspired by the interesting paper due to Savo \cite{Savo14}. The key point  
is to establish a sharp inequality for the Bochner operator in 
terms of the second fundamental form of any immersed submanifold in a Riemannian manifold whose 
curvature operator is bounded from below by a constant. Moreover, we determine the structure of the second 
fundamental form at points where equality holds. All the above mentioned, combined with Bochner technique 
and Morse theory, enable us to provide homology vanishing results, bounds for the Betti numbers, or determine 
the submanifolds that satisfy the pinching condition \eqref{1} regardless codimension. Examples 
show that the pinching  condition is sharp.

\section{The results}
This section is devoted to the statements of our results. Throughout the paper, all submanifolds under consideration are assumed to be oriented. Our first global result is stated as follows.

\begin{theorem}\label{mainthm1}
Let $M^n, n\geq 3$, be a compact $n$-dimensional Riemannian manifold isometrically immersed in a 
Riemannian manifold whose curvature operator is bounded from below by a constant $c$. Suppose that 
inequality \eqref{1} is satisfied for an integer $1\leq p\leq n/2$ and $H^2+c\geq 0$. 
The following assertions hold:
\begin{enumerate}[(i)]
\item The Betti numbers of $M^n$ satisfy $b_i(M^n)\leq \binom n i$ for each $p\leq i\leq n-p$.
\item If strict inequality holds in \eqref{1} at a point, then $b_i(M^n)=0$ for all $p\leq i\leq n-p$. 
If in addition $p=1$, then $M^n$ is a real homology sphere, i.e., $H^i(M^n;\R)=0$ for each $1\leq i\leq n-1$.
\item If $b_p(M^n)>0$, then equality holds in \eqref{1} everywhere on $M^n$ and every harmonic 
$p$-form is parallel. In particular, $M^n$ supports a nontrivial parallel $p$-form. 
\end{enumerate}
\end{theorem}

The above extends previous results due to Cui and Sun \cite{CS19} and Savo \cite[Th. 8]{Savo14}. 

Our second main result is about (not necessarily compact) submanifolds carrying a nontrivial parallel 
$p$-form and can be formulated as follows (see \cite{Gr} for related results).

\begin{theorem}\label{thm2}
Let $M^n, n\geq 3$, be a Riemannian manifold isometrically immersed in 
a  Riemannian manifold whose curvature operator is bounded from below by a constant $c$ such that 
$H^2+c\geq 0$. If $M^n$ supports a nontrivial harmonic $p$-form of constant length (in particular, 
a parallel $p$-form), then $S\geq a(n,p,H,c).$
\end{theorem}

Since K\"ahler manifolds support non-trivial parallel forms in all even degrees, we obtain the following:

\begin{corollary}\label{kaehler}
If a K\"ahler manifold $M^{2m}$ is isometrically immersed in a Riemannian manifold whose
curvature operator is bounded from below by a constant $c$ such that $H^2+c\geq 0$, then
$S\geq a\big(2m,2[m/2],H,c\big).$
\end{corollary}

\medskip
In the above results, the assumption that the curvature operator of the ambient manifold is bounded from 
below can be weakened by merely assuming that the induced Bochner operator is bounded from below.

Now we focus on submanifolds that satisfy the pinching condition \eqref{1} in complete simply connected space 
forms $\Q^m_c$ of constant sectional curvature $c$, that is, the Euclidean 
space $\R^m$, the sphere $\Sf^m$ or the hyperbolic space $\mathbb H^m$, according to whether $c=0, c=1$  or $c =-1$, 
respectively. The next result deals with minimal submanifolds in spheres and extends the results in \cite{an, leung83}. 

\begin{theorem}\label{thm.b}
Let $f\colon M^n\to \Sf^m, n\geq 3$, be a minimal isometric immersion of a compact Riemannian manifold 
$M^n$. If $S\leq n$ on $M^n$, then one of the following assertions holds:
\begin{enumerate}[(i)]
\item $M^n$ is a real homology sphere, with finite fundamental group, that admits a Riemannian metric of positive 
Ricci curvature. 
\item $f(M^n)$ is isometric to the Clifford torus $\mathbb T^n_p(\sqrt{p/n})$ 
for an integer $1\leq p\leq n/2$. 
\item The immersion $f$ is the standard embedding $\psi\colon \CP^2_{4/3}\to \Sf^7$ of the complex 
projective plane of constant holomorphic curvature $4/3$.
\end{enumerate}
\end{theorem}

In particular, for minimal submanifolds of dimension three and four, we have the following result that improves 
the bound in \cite[Th. 3]{ybshen}:

\begin{corollary}\label{cor}
Let $f\colon M^n\to \Sf^m$ be a minimal isometric immersion of a compact Riemannian manifold $M^n$ with 
$S\leq n$ on $M^n$. 
\begin{enumerate}[(i)]
\item If $n=3$, then either  $f(M^3)$ is isometric to the Clifford torus $\mathbb T^3_1(\sqrt{1/3})$, 
or $M^3$ is diffeomorphic to a spherical space form.
\item If $n=4$, then $f(M^4)$ is isometric to the Clifford torus $\mathbb T^4_p(\sqrt{p/n})$, $p=1,2$, 
the immersion $f$ is the standard embedding $\psi\colon \CP^2_{4/3}\to \Sf^7$, or the universal cover of 
$M^4$ is homeomorphic to $\Sf^4$. 
\end{enumerate}
\end{corollary}

We recall that the Cartan minimal hypersurface in $\mathbb{S}^4$ satisfies $S=6$, while  the minimal submanifold in $\mathbb{S}^5$, constructed in \cite{dv} as a circle bundle over a flat minimal torus, has $S=8$.  In addition, the minimal immersion of the $n$-dimensional complex projective $\CP^2_{2n/(n+1)}$ space of constant holomorphic curvature $2n/(n+1), n>2,$ into $\mathbb{S}^{n(n+2)}$ (see \cite{W}) satisfies $S=2n(n-1)$. These examples show that Theorem \ref{thm.b}  and Corollary \ref{cor} are optimal.

Our results for not necessarily minimal submanifolds of spheres can be stated as follows:

\begin{theorem}\label{thm.a}
Let $f\colon M^n\to \Sf^m, n\geq 3$, be an isometric immersion of a compact 
Riemannian manifold $M^n$. Suppose that inequality \eqref{1} is satisfied for an integer $1\leq p\leq n/2$. 
Then either $b_i(M^n)=0$ for every $p\leq i\leq n-p$, or the following assertions hold:
\begin{enumerate}[(i)]
\item If $b_p(M^n)>0$, then $M^n$ has the homotopy type of a CW-complex 
with cells only in dimensions $0,p,n-p$ or $n$. In particular, the homology groups of $M^n$ satisfy 
$H_i (M^n;G)=0$ for each $i\neq 0,p,n-p,n$, where $G$ is any coefficient group.
\item If $b_q(M^n)>0$ for an integer $p<q\leq n/2$, then $f$ is minimal and either 
$f(M^n)$ is isometric to the Clifford torus $\mathbb T^n_q(\sqrt{q/n})$, or $f$ is the standard embedding 
as in Theorem \ref{thm.b}(iii).
\end{enumerate}
If in addition $H>0$ on $M^n$, then the following facts hold:
\begin{enumerate}[(i)]
\item[(iii)] If $b_p(M^n)>0$ and $p<n/2$,  
then $f(M^n)$ is isometric to a torus $\mathbb T^n_p(r)$ with $r>\sqrt{p/n}$. 
\item[(iv)] If $p<(n-1)/2$, then either $H_i (M^n;\Z)=0$ for all $p\leq i< n-p$, or $f(M^n)$ is as in (iii) above. 
\end{enumerate}
\end{theorem}

The above result extends the sphere theorem proved in \cite{sx97.2} and \cite[Th. 9]{Savo14}. The standard 
embedding $\psi\colon \CP^2_{4/3}\to \Sf^7$ justifies the necessity of the assumptions in parts (iii) and (iv). 

\begin{theorem}\label{cor0}
Let $f\colon M^n\to \Sf^m, n\geq 3,$ be an isometric immersion of a compact Riemannian manifold $M^n$. If inequality \eqref{1} 
holds for $p=1$, then one of the following assertions holds:
\begin{enumerate}[(i)]
\item $M^n$ is a real homology sphere, with finite fundamental group, that admits a Riemannian metric of positive 
Ricci curvature. Furthermore, $M^n$ is diffeomorphic 
to a spherical space form if $n=3$ and its universal cover is homeomorphic to $\Sf^4$ if $n=4$.
\item $M^n$ has the homotopy type of a CW-complex with cells only in dimensions $0,1,n-1,n$.
\item $f(M^n)$ is isometric to the Clifford torus $\mathbb T^n_q(\sqrt{q/n})$ for an integer $1<q\leq n/2$, or $f$ 
is the standard embedding as in Theorem \ref{thm.b}(iii).
\end{enumerate}
If in addition $H>0$ on $M^n$, then either
$M^n$ is homeomorphic to $\Sf^n$, or $f(M^n)$ is isometric to a torus $\mathbb T^n_1(r)$ with $r>1/\sqrt{n}$.
\end{theorem}

The pinching condition \eqref{1}  is sharp in both Theorems \ref{thm.a} and \ref{cor0}. In fact, let $g$ be  the Cartan minimal hypersurface in $\mathbb{S}^4$, the minimal submanifold in $\mathbb{S}^5$ constructed in \cite{dv}, or the minimal immersion of the $n$-dimensional complex projective $\CP^2_{2n/(n+1)}$ space of constant holomorphic curvature $2n/(n+1), n>2,$ into $\mathbb{S}^{n(n+2)}$. For each $0<\rho<1$ we consider the isometric immersion $f_\rho=j\circ G_\rho \circ g$, where $G_\rho$ is the homothety of center at the origin and ratio $\rho$, and  $j$ is the standard inclusion of the sphere of radius $\rho$  as a hypersurface in the unit sphere. It is easy to see that for $\rho$ sufficiently close to $1$ the nonminimal  immersion $f_\rho$ satisfies the strict reverse inequality in \eqref{1}.

The following corollary improves a result due to Lawson and Simons \cite[Cor. 2]{LS73}. 

\begin{corollary}\label{cor1}
Let $f\colon M^n\to \Sf^m, n\geq 3,$ be an isometric immersion of a compact Riemannian manifold $M^n$. If 
$S\leq 2\sqrt{n-1}$ on $M^n$, then either $f(M^n)$ is isometric to a torus $\mathbb T^n_1(r)$ with 
appropriate $r>1/\sqrt{n}$, or $M^n$ is a real homology sphere of positive Ricci curvature. Furthermore, 
$M^n$ is diffeomorphic to a spherical space form if $n=3$ and its universal cover is homeomorphic to $\Sf^4$ 
if $n=4$.
\end{corollary}

For Euclidean and spherical submanifolds, that are contained in an open hemisphere, we prove the following result:

\begin{theorem}\label{thm.A1}
Let $f\colon M^n\to \Q_c^m, n\geq 3$, be an isometric immersion of a compact Riemannian manifold $M^n$ 
and let $c=0,1$. If $c=1$, suppose further that $f(M^n)$ is contained in an open hemisphere. If inequality \eqref{1} holds 
for an integer $1\leq p\leq n/2$, then the following assertions hold:
\begin{enumerate}[(i)]
\item Either $b_i(M^n)=0$ for each $p\leq i\leq n-p$, or $p=n/2$ 
and $M^n$ has the homotopy type of a CW-complex with cells only in dimensions $0,n/2$ or $n$. 
\item If $p=1$, then $M^n$ is a real homology sphere, with finite fundamental group, that admits a 
Riemannian metric of positive Ricci curvature. Furthermore, $M^n$ is diffeomorphic 
to a spherical space form if $n=3$ and its universal cover is homeomorphic to $\Sf^4$ if $n=4$.
\item If $H>0$ on $M^n$ and $p<(n-1)/2$ when $c=0$, then $H_i (M^n;\Z)=0$ for each $p\leq i < n-p$.
\end{enumerate}
\end{theorem}

The assumption in the above result that the submanifold is contained in an open hemisphere is necessary, as 
shown by the Clifford tori in Theorem \ref{thm.b}(ii). 
It is worth noticing that the even dimensional torus $\mathbb T^n_{n/2}(1/\sqrt{2})$, viewed as a 
submanifold of the Euclidean space, satisfies equality in the pinching condition \eqref{1}. 
The above result extends Theorem 1.6 in \cite{CS19} and Theorem 1.2 in \cite{lin15}. 
A stronger result was proved in \cite{guxu12} for the case $p=1$ under the additional assumption that the mean 
curvature is positive everywhere.

The following improves the bound given in \cite{fuxu08} for submanifolds of a hyperbolic space.

\begin{theorem}\label{thm.d}
Let $M^n, n\geq 3$, be a compact Riemannian manifold isometrically immersed in a hyperbolic 
space $\Hy^m$ with mean curvature $H\geq 1$. If inequality \eqref{1} holds 
for an integer $1\leq p< n/2$, then $b_i(M^n)=0$ for every $p\leq i\leq n-p$. 
\end{theorem}

Next, we provide an integral bound for the Bochner operator in terms of the mean 
curvature, the length of the traceless part $\Phi$ of the second fundamental form 
and the Betti numbers $b_{i}(M^n;\mathbb{F})$ over an arbitrary coefficient field $\mathbb F$.

\begin{theorem}\label{thm.integral}
Given integers $n\geq 4,\, 1\leq p< n/2$ and $k\geq 1$, there exists a positive constant $c(n,p,k)$ 
such that if $M^n$ is a compact Riemannian manifold isometrically immersed 
in $\R^{n+k}$, then the lowest eigenvalue $\varrho_p$ of its Bochner operator $\B^{[p]}$ satisfies
$$
\int_{M^n} \big| \varrho_p -F_{n,p}(H, \|\Phi\|) \big|^{n/2} dM
\geq c(n,p,k)\sum_{p<i<n-p} b_{i}(M^n;\mathbb{F})
$$
for any coefficient field $\mathbb{F}$, where $F_{n,p}$ is the function given by 
$$
F_{n,p}(x, y)=\frac{p(n-p)}{n} \big(nx^2 -\frac{n(n-2p)}{\sqrt{np(n-p)}}x y-y^2\big).
$$ 
Moreover, if 
$$
\int_{M^n} \big| \varrho_p -F_{n,p}(H, \|\Phi\|) \big|^{n/2} dM< c(n,p,k),
$$ 
then $M^n$ has the homotopy type of a CW-complex with no cells of dimension 
$p<i<n-p$. In particular, if $p=1$, then the fundamental group $\pi_1(M^n)$ is a free group on 
$b_1(M^n)$ generators, and $M^n$ is homeomorphic to $\Sf^n$ if $\pi_1(M^n)$ is finite.
\end{theorem}

\section{Algebraic preliminaries and auxiliary results}

This section is devoted to some algebraic preliminaries inspired by the ideas developed by Savo \cite{Savo14}.
Let $V$ be a real $n$-dimensional vector space, $n\geq 3$, equipped with a positive definite 
inner product denoted by $\langle\cdot,\cdot\rangle$. We denote by $\mathrm{End}(V)$ 
the set of all self-adjoint endomorphisms of $V$. 

For every integer $1\leq p\leq n$, we consider the set $I_p$ of $p$-multi-indices 
$$
I_p=\left\{\{i_1,\dots,i_p\} : 1\leq i_1<\cdots<i_p\leq n\right\}.
$$
Let $A\in \mathrm{End}(V)$ be an arbitrary endomorphism of $V$ with eigenvalues $k_1,\dots,k_n$. 
For each $a=\{i_1,\dots,i_p\}\in I_p$, the associated $p$-{\it algebraic curvature} $K_a$ of $A$ is
the number given by
$$
K_a= k_{i_1}+\cdots+k_{i_p}.
$$

Let $\Lambda^p V^*$ be the $\binom np$-dimensional real vector space, defined as the
$p$-th exterior power of the dual vector space $V^*=\Hom(V,\R)$ of $V$. 
To each $A\in \End(V)$, we associate the endomorphism $T_A^{[p]}\in\End(\Lambda^p V^*)$ defined by 
$$
T_A^{[p]}=(\tr A)A^{[p]}-A^{[p]}\circ A^{[p]},
$$
where $A^{[p]}\in \End(\Lambda^p V^*)$ is given by 
$$
A^{[p]}\omega(v_1,\dots,v_p)=\sum_{i=1}^p \omega(v_1,\dots,Av_i,\dots,v_p),
$$
with $\omega\in \Lambda^p V^*$ and $v_1,\dots,v_p\in V$. The endomorphism $T_A^{[p]}$ is self-adjoint 
with respect to the natural inner product $\<\cdot, \cdot \>$ in $\Lambda^p V^*$.
It follows from \cite[Lem. 2]{Savo14} that the $\binom np$ eigenvalues $\lambda_a(T_A^{[p]})$ of $\Tp$ 
are given by 
\begin{equation}\label{lambda}
\lambda_a(T_A^{[p]})=K_a K_{\star a},\; \; a=\{i_1,\dots,i_p\}\in I_p,
\end{equation}
where $\star a \in I_p$ is defined by $\star a=\{1,\dots,n\}\smallsetminus\{i_1,\dots,i_p\}$.

The lowest eigenvalue of $T_A^{[p]}$ is given by
$$
\min_{a\in I_p} \lambda_a(T_A^{[p]})
=\mathop{\min_{\omega\in \Lambda^p V^*}}_{\|\omega\|=1}\<T_A^{[p]}\omega,\omega\>.
$$

We denote by $\mathring A=A-H_A {\rm id}_V$ the traceless part of every $A\in \End(V)$, where 
$H_A=(1/n)\, \tr A$ and ${\rm id}_V\colon V\to V$ being the identity map on $V$.

The inequality in the following lemma was essentially proved in \cite[Lem. 14]{Savo14}.

\begin{lemma}\label{lem1}
For every $A\in \mathrm {End}(V)$ with $\tr A\geq 0$ and each integer $1\leq p\leq n/2$, the lowest eigenvalue 
of $T_A^{[p]}\in\mathrm {End}(\Lambda^p V^*)$ satisfies
\be\label{basineq1}
\min_{a\in I_p} \lambda_a(T_A^{[p]})\geq \frac{1}{4}(\tr A)^2-\frac{1}{4}\Big(\frac{n-2p}{n}\, \tr A 
+\sqrt{4p(n-p)/n}\, \|\mathring A\|\Big)^2.
\ee
If equality holds, then $A$ has at most two distinct eigenvalues 
$\lambda$ and $\mu$ with multiplicities $p$ and $n-p$, respectively. If in addition $\mu=0$ 
and $p<n/2$, then $A=0$.
\end{lemma}

\proof
Let $k_1,\dots,k_n$ be the eigenvalues of $A$. We distinguish two cases. 

\medskip
{\it Case 1}. Suppose that $\tr A=0$ and let $a=\{i_1,\dots,i_p\} \in I_p$. It follows from \eqref{lambda} 
that the eigenvalue $\lambda_a(T_A^{[p]})$ of $T_A^{[p]}$ satisfies 
\bea
\lambda_a(T_A^{[p]})&=& -\big(\sum_{i\in a} k_i\big)^2\geq -p\sum_{i\in a} k_i^2, \\
\lambda_a(T_A^{[p]}) &=& -\big(\sum_{i\in \star a} k_i\big)^2\geq -(n-p)\sum_{i\in \star a} k_i^2.
\eea
Summing up, we obtain
\be\label{ineq1}
\lambda_a(T_A^{[p]})\geq -\frac{p(n-p)}{n} \|A\|^2
\ee
which proves the desired inequality in this case. If equality holds in \eqref{basineq1}, then all the 
above inequalities become equalities, and this implies that $A$ has at most two distinct eigenvalues
$$
\lambda=k_{i_1}=\cdots=k_{i_p}\; \; \text{and}\; \; \mu=k_{i_{p+1}}=\cdots=k_{i_{n}}
$$
with multiplicities $p$ and $n-p$, respectively, where $\star a=\{i_{p+1},\dots,i_{n}\}$.

\medskip
{\it Case 2}. Suppose that $\tr A>0$. The eigenvalues of the traceless part $\mathring A$ of $A$ are 
$k_1-(\tr A)/n,\dots,k_n-(\tr A)/n$. It follows from \eqref{lambda}
 that, for each $a=\{i_1,\dots,i_p\} \in I_p$, the eigenvalue $\lambda_a(T_{\mathring A}^{[p]})$ of 
$T_{\mathring A}^{[p]}$ is given by
\bea
\lambda_a(T_{\mathring A}^{[p]}) \!\!\!&=&\!\!\! \big(K_a-\frac{p}{n}\, \tr A\big)\big(K_{\star a} -\frac{n-p}{n}\, \tr A\big) \nonumber \\
\!\!\!&=&\!\!\! \lambda_a(T_{A}^{[p]})+\frac{p(n-p)}{n^2}\, (\tr A)^2-\frac{n-p}{n}K_a\, \tr A-\frac{p}{n}K_{\star a}\, \tr A. \nonumber
\eea
Using that $K_a+K_{\star a}=\tr A$, we obtain 
\begin{eqnarray*}
\lambda_a(T_{\mathring A}^{[p]}) &=& \lambda_a(T_A^{[p]})-\frac{p^2}{n^2}\left(\tr A\right)^2+\frac{2p-n}{n}K_a\, \tr A, \\
\lambda_a(T_{\mathring A}^{[p]})&=&\lambda_a(T_A^{[p]})-\frac{(n-p)^2}{n^2}\left(\tr A\right)^2+\frac{n-2p}{n}K_{\star a}\, \tr A. 
\end{eqnarray*}
Summing up, we have
$$
2\lambda_a(T_{\mathring A}^{[p]})=2\lambda_a(T_A^{[p]})-\frac{(n-p)^2+p^2}{n^2}\left(\tr A\right)^2+\frac{n-2p}{n} \left(K_{\star a}-K_a\right) \tr A.
$$
Since $\mathring A$ is trace free, inequality \eqref{ineq1} applies and yields 
$$
\lambda_a(T_{\mathring A}^{[p]})\geq -\frac{p(n-p)}{n} \|\mathring A\|^2.
$$
Combining the above, we obtain
\bea
2\lambda_a(T_A^{[p]})-\frac{(n-p)^2+p^2}{n^2}\left(\tr A\right)^2+\frac{2}{n}p(n-p) \|\mathring A\|^2\\ 
\geq -\frac{1}{n} (n-2p) \left|K_{\star a}-K_a\right|\tr A. 
\eea
Using \eqref{lambda}, it follows that
$$
|K_{\star a}-K_a|=\big((\tr A)^2-4\lambda_a(T_A^{[p]})\big)^{1/2},
$$
and consequently 
\bea
2\lambda_a(T_A^{[p]})-\frac{(n-p)^2+p^2}{n^2}(\tr A)^2+\frac{2}{n}p(n-p) \|\mathring A\|^2\\ \nonumber 
\geq -\frac{1}{n}\left(n-2p\right) \tr A \big((\tr A)^2-4\lambda_a(T_A^{[p]})\big)^{1/2}. \nonumber 
\eea
Setting $\rho=\big((\tr A)^2-4\lambda_a(T_A^{[p]})\big)^{1/2}$,  
the above is written as 
$$
n^2\rho^2-2n(n-2p)\, \tr A\, \rho+(n-2p)^2(\tr A)^2-4np(n-p) \|\mathring A\|^2\leq 0.
$$
Thus
$$
\rho\leq \frac{1}{n}\left(n-2p\right) \tr A +\sqrt{4p(n-p)/n}\, \|\mathring A\|,
$$
or equivalently
$$
\lambda_a(T_A^{[p]})\geq \frac{1}{4}(\tr A)^2-\frac{1}{4} \big(\frac{n-2p}{n}\, \tr A +\sqrt{4p(n-p)/n}\,\|\mathring A\|\big)^2.
$$
Inequality \eqref{basineq1} follows immediately from the above. If equality holds in \eqref{basineq1}, then 
all the above inequalities become equalities, and this implies that $A$ has at most two distinct eigenvalues
$$
\lambda=k_{i_1}=\cdots=k_{i_p}\; \; \text{and}\; \; \mu=k_{i_{p+1}}=\cdots=k_{i_{n}}
$$
with multiplicities $p$ and $n-p$, respectively, where $\star a=\{i_{p+1},\dots,i_{n}\}$.

\medskip
Next we suppose that equality holds in \eqref{basineq1} and 
$$
\lambda_a(T_A^{[p]})=\min_{b \in I_p} \lambda_b(T_A^{[p]})
$$
for some $a=\{i_1,\dots,i_p\}\in I_p$. In both Cases 1 and 2, the endomorphism $A$ has at most two 
distinct eigenvalues
$$
\lambda=k_{i_1}=\cdots=k_{i_p}\; \; \text{and}\; \; \mu=k_{i_{p+1}}=\cdots=k_{i_n}
$$
with multiplicities $p$ and $n-p$, respectively, where $\star a=\{i_{p+1},\dots,i_n\}$. A direct computation 
shows that the eigenvalues of $T_A^{[p]}$ are given by
$$
\lambda_b(T_A^{[p]})=\left((p-\mathsf n_b)\lambda+\mathsf n_b\mu\right)\left(\mathsf n_b\lambda+(n-p-\mathsf n_b)\mu \right) 
$$
for each $b=\{j_1,\dots,j_p\} \in I_p$ with $1\leq j_1<\dots < j_p\leq n$. Here $0\leq \mathsf n_b\leq p$ 
is the number of common elements of $b$ and $a=\{i_1,\dots,i_p\}$.

In addition, we suppose that $\mu=0$ and $p<n/2$. From the above it follows that
$$
\min_{b\in I_p} \lambda_b(T_A^{[p]})=0.
$$
Computing the right hand side in \eqref{basineq1} in terms of the eigenvalues of $A$, 
and since equality holds in \eqref{basineq1}, we obtain
$$
\min_{b\in I_p} \lambda_b(T_A^{[p]})=\frac{2p^2}{n^2}(n-p)(2p-n)\lambda^2.
$$
It is now obvious that $A=0$.
\qed

\medskip
Let $W$ be a $k$-dimensional real vector space equipped with a positive definite 
inner product which, by abuse of notation, is again denoted by $\langle\cdot,\cdot\rangle$.
We denote by $\Hom(V\times V,W)$ the space of all bilinear forms with values in $W$ and by 
$\mathrm{Sym}(V\times V,W)$ the subspace that consists of all symmetric bilinear forms. 
The space $\mathrm{Sym}(V\times V,W)$ can be viewed as a complete metric space with respect 
to the usual Euclidean norm $\Vert\cdot\Vert$. For each $\beta\in \mathrm{Sym}(V\times V,W)$, 
we define the map
$$
\beta^\sharp\colon W\rightarrow \mathrm{End}(V),\; \; \xi\mapsto \beta^\sharp(\xi)
$$
such that
$$
\langle \beta^\sharp(\xi)v_1,v_2\rangle=\langle \beta(v_1,v_2),\xi\rangle\;\;\text{for all}\;\; v_1,v_2\in V.
$$

For every integer $1\leq p \leq n$, we define the map 
$$
\mathsf B^{[p]} \colon \mathrm{Sym}(V\times V,W) \to \End(\Lambda^p V^*)
$$
given by
$$
\mathsf B^{[p]}(\beta)=\sum_{i=1}^k T_{\beta^\sharp(\xi_i)}^{[p]},
$$
where $\xi_1,\dots,\xi_k$ is an arbitrary orthonormal basis of $W$. Observe that $\mathsf B^{[p]}(\beta)$ 
is a self-adjoint endomorphism of $\Lambda^p V^*$. Its lowest eigenvalue is denoted by $\varrho_p(\beta)$.

For convenience, we set $\mathsf H_\beta= \left(1/n\right) \mathrm{tr} \beta$ and denote by 
$\mathring{\beta}= \beta-\<\cdot, \cdot \>\mathsf H_\beta$ 
the traceless part of each $\beta\in \mathrm{Sym}(V\times V,W)$. The image of $\beta$ is 
the subspace $\mathrm{Im}\, \beta \subset W$ given by
$$
\mathrm{Im}\, \beta=\span\left\{\beta (v_1,v_2): v_1,v_2\in V\right\}.
$$ 

The following will be an important tool in our proofs.

\begin{proposition}\label{lem2}
For every $\beta\in \mathrm{Sym}(V\times V,W)$ and each integer $1\leq p\leq n/2$,
the lowest eigenvalue $\varrho_p(\beta)$ of $\mathsf B^{[p]}(\beta)\in\mathrm{End}\left(\Lambda^p V^*\right)$ 
satisfies
\be\label{ineq2}
\varrho_p(\beta)
\geq \frac{p(n-p)}{n} \big(n\|\mathsf H_{\beta}\|^2 -\frac{n(n-2p)}{\sqrt{np(n-p)}}\, \|\mathsf H_{\beta}\|\|\mathring{\beta}\|-\|\mathring{\beta}\|^2\big).
\ee
Moreover, if equality holds in \eqref{ineq2} for some $\beta$, then the following hold:

\begin{enumerate}[(i)]
\item For every unit vector $u\in W$, $\beta^\sharp(u)$ has at most two distinct 
eigenvalues of multiplicities $p$ and $n-p$. In addition, if $p<n/2$ and the eigenvalue with 
multiplicity $n-p$ vanishes, then $\beta^\sharp(u)=0$.

\item If $\mathsf H_\beta\neq 0$ and $p<n/2$, then $\mathrm{Im}\, \beta=\span\left\{\mathsf H_\beta\right\}$.
\end{enumerate}
\end{proposition}

\proof
Let $\beta\in \mathrm{Sym}(V\times V,W)$. If $\dim W=1$, then the proof follows directly from 
Lemma \ref{lem1}. We suppose that $\dim W>1$ and distinguish two cases. 

\smallskip
{\it Case 1}. Suppose that $\mathsf H_\beta\neq0$ and consider an orthonormal basis 
$\xi_1,\dots, \xi_k$ of $W$ such that $\xi_1=\mathsf H_\beta/\|\mathsf H_\beta\|$.
Then we have $\tr \beta^\sharp(\xi_i)=0$ for each $2\leq i\leq k$.
From Lemma \ref{lem1}, it follows that the eigenvalues of $T_{\beta^\sharp(\xi_i)}^{[p]},1\leq i\leq k$, 
satisfy
\be\label{ineq5}
\min_{a\in I_p}\lambda_a(T_{\beta^\sharp(\xi_1)}^{[p]}) \geq
\frac{n^2}{4} \|\mathsf H_\beta\|^2
-\frac{1}{4}\big((n-2p)\|\mathsf H_\beta\|
+\sqrt{4p(n-p)/n}\, \|\mathring{\beta}^\sharp(\xi_1)\|\big)^2
\ee
and 
\be\label{ineq6}
\min_{a\in I_p}\lambda_a(T_{\beta^\sharp(\xi_i)}^{[p]}) \geq
-\frac{p(n-p)}{n}\,\|\mathring{\beta}^\sharp(\xi_i)\|^2\;\; \text{for each}\;\; 2\leq i\leq k.
\ee
Note that 
$$
\|\mathring{\beta}^\sharp(\xi_1)\|^2 = \|\beta^\sharp(\xi_1)\|^2-n\|\mathsf H_\beta\|^2 \; \; \text{and}\; \;
\mathring{\beta}^\sharp(\xi_i) =\beta^\sharp(\xi_i) \;\; \text{for all}\;\; 2\leq i\leq k. 
$$
Summing up \eqref{ineq5} and \eqref{ineq6} and since 
$$
\varrho_p(\beta)=\mathop{\min_{\omega\in \Lambda^p V^*}}_{\|\omega\|=1}\<\mathsf B^{[p]}(\beta)\omega,\omega\>=\min_{a\in I_p}\lambda_a(\mathsf B^{[p]}(\beta)),
$$
we obtain 
\bea
\varrho_p(\beta)\!\!\!&\geq&\!\!\!
p(n-p)\|\mathsf H_\beta\|^2 
-\frac{p(n-p)}{n}\big(\|\beta\|^2-n\|\mathsf H_\beta\|^2\big) \nonumber \\
\!\!\!& &\!\!\!-(n-2p)\|\mathsf H_\beta\|\sqrt{p(n-p)/n}\, \sqrt{\|\beta^\sharp(\xi_1)\|^2-n\|\mathsf H_\beta\|^2} \nonumber\\
\!\!\!&\geq&\!\!\! p(n-p)\|\mathsf H_\beta\|^2-\frac{p(n-p)}{n}\big(\|\beta\|^2-n\|\mathsf H_\beta\|^2\big) \nonumber \\
\!\!\!& &\!\!\!-(n-2p)\|\mathsf H_\beta\|\sqrt{p(n-p)/n}\, \sqrt{\|\beta\|^2-n\|\mathsf H_\beta\|^2}, \nonumber 
\eea
and this is inequality \eqref{ineq2}. 

Now, we suppose that equality holds in \eqref{ineq2} for some $\beta\in \mathrm{Sym}(V\times V,W)$. 
Since all inequalities above become equalities, we obtain $\beta^\sharp(\xi_i)=0$ for each $2\leq i\leq k$ 
provided that $p<n/2$. Thus $\mathrm{Im}\, \beta=\span\left\{\mathsf H_\beta\right\}$ if $p<n/2$, and the 
rest of the proof for this case follows from Lemma \ref{lem1}.

\smallskip
{\it Case 2}. Suppose that $\beta\in \mathrm{Sym}(V\times V,W)$ is traceless. Let $u\in W$ be an arbitrary 
unit vector. We choose an orthonormal base $\xi_1,\dots,\xi_k$ of $W$ such that $\xi_1=u$. Since 
$\tr \beta^\sharp(\xi_i)=0$ for every $1\leq i\leq k$, from Lemma \ref{lem1} it follows that 
$$
\min_{a\in I_p}\lambda_a(T_{\beta^\sharp(\xi_i)}^{[p]})\geq 
-\frac{p(n-p)}{n}\, \|\beta^\sharp(\xi_i)\|^2\;\; \text{for each}\;\; 1\leq i\leq k.
$$
Summing over $i$, we obtain 
$$
\varrho_p(\beta)=\min_{a\in I_p}\lambda_a(\mathsf B^{[p]}(\beta))
\geq -\frac{p(n-p)}{n}\sum_{i=1}^k\|\beta^\sharp(\xi_i)\|^2 = -\frac{p(n-p)}{n}\, \|\beta\|^2
$$
and this proves inequality \eqref{ineq2}.

If equality holds in \eqref{ineq2} for some $\beta\in \mathrm{Sym}(V\times V,W)$, 
then the above inequalities become equalities and the proof follows from Lemma \ref{lem1}.
\qed

\medskip
Next we assume that the dimension of $V$ is $n\geq 4$. For each integer $1\leq p< n/2$ and for every 
$\beta\in \mathrm{Sym}(V\times V,W)$, we denote by $\Lambda_p(\beta)$ the 
subset of the unit $(k-1)$-sphere $\Sf^{k-1}$ in $W$ given by
$$
 \Lambda_p(\beta)=\left\{u\in\mathbb{S}^{k-1}:p< \mathrm{Index}\, \beta^\sharp(u)<n-p\right\}.
$$

The inequality that follows is crucial for the proof of Theorem \ref{thm.integral}.

\begin{proposition}\label{int.lemma}
Given integers $n\geq 4,\, k\geq 1$ and $1\leq p< n/2$, there exists a constant 
$\varepsilon(n,k,p)>0$ such that the following inequality holds
\bea
\varrho_p(\beta)
- \frac{p(n-p)}{n} \big(n\|\mathsf H_{\beta}\|^2 -\frac{n(n-2p)}{\sqrt{np(n-p)}}\, \|\mathsf H_{\beta}\|\|\mathring{\beta}\|
-\|\mathring{\beta}\|^2\big)\\ \nonumber
\geq 
\varepsilon(n,k,p)\Big(\int_{\Lambda_p(\beta)}\vert\det\beta^\sharp (u)\vert dS_u\Big)^{2/n}
\eea
for all $\beta\in \mathrm{Sym}(V\times V,W)$, where $dS_u$ stands for the volume element 
of the unit sphere in $W$.
\end{proposition}

\proof
We consider the functions 
$\phi_p,\psi_p\colon \mathrm{Sym}(V\times V,W)\rightarrow \mathbb{R}$ defined by
$$
\phi_p(\beta) = 
\varrho_p(\beta)
- \frac{p(n-p)}{n} \big(n\|\mathsf H_{\beta}\|^2 -\frac{n(n-2p)}{\sqrt{np(n-p)}}\, \|\mathsf H_{\beta}\|\|\mathring{\beta}\|
-\|\mathring{\beta}\|^2 \big)
$$
and 
$$
\psi_p(\beta)=\int_{\Lambda_p(\beta)}\vert\det\beta^\sharp (u)\vert dS_u.
$$
From Proposition \ref{lem2} we know that the function $\phi_p$ is nonnegative. To prove the 
desired inequality, it is sufficient to show that $\phi_p$ attaints a positive minimum on the level set
$$
\Sigma_{n,k,p}=\left\{\beta\in\mathrm{Sym}(V\times V,W):\psi_p(\beta)=1\right\}.
$$
Let $\{\beta_m\}$ be a sequence in $\Sigma_{n,k,p}$ such that 
$$\lim_{m\rightarrow\infty}\phi_p(\beta_m)=\inf \phi_p(\Sigma_{n,k,p})\geq 0.$$

\smallskip
{\it Claim I: The sequence $\{\beta_m\}$ is bounded.} Suppose to the contrary that 
there exists a subsequence of $\{\beta_m\}$, which by abuse of notation is again 
denoted by $\{\beta_m\}$, such that $\lim_{m\rightarrow \infty}\Vert\beta_m\Vert=\infty$. 
Since $\beta_m\neq0$ for all $m\in\N$, we may write $\beta_m=\Vert\beta_m\Vert\hat{\beta}_m$, where 
$\Vert\hat{\beta}_m\Vert=1$, and assume that $\{\hat{\beta}_m\}$ converges to some 
$\hat{\beta}\in\mathrm{Sym}(V\times V,W)$ with $\Vert\hat{\beta}\Vert=1$. It can be easily checked that 
$\phi_p$ is homogeneous of degree $2$, that is 
$$
\phi_p(t\beta)=t^2\phi_p(\beta)\; \; \text{for all}\;\; t>0.
$$
Therefore, we obtain $\phi_p(\hat{\beta}_m)=\phi_p(\beta_m)/\Vert\beta_m\Vert^2$. 
Thus $\lim_{m\rightarrow\infty}\phi_p(\hat{\beta}_m)=0$ and consequently 
$\phi_p(\hat{\beta})=0$, which means that $\hat{\beta}$ satisfies equality in inequality \eqref{ineq2} 
in Proposition \ref{lem2}. To reach a contradiction, we distinguish two cases.

\smallskip
{\it Case 1}. Suppose that $\mathsf H_{\hat{\beta}}\neq 0$. Proposition \ref{lem2} implies that 
$\hat{\beta}(\cdot , \cdot)=\< {\hat\beta}^\sharp(\hat\xi) \cdot , \cdot \>\, \hat\xi$ 
and ${\hat\beta}^\sharp(\hat\xi)$ has at most two distinct eigenvalues $\hat\lambda$ and $\hat \mu$ with 
multiplicities $p$ and $n-p$, respectively, where 
$\hat \xi=\mathsf H_{\hat{\beta}}/\|\mathsf H_{\hat{\beta}}\|$. 
Clearly $\hat\mu\neq0$, since otherwise $\hat\beta=0$, and this is a contradiction. 

Since $\beta_m\in\Sigma_{n,k,p}$, there exists an open subset $\hat{\CU}_m\subset \Sf^{k-1}\subset W$ 
such that $\hat{\CU}_m\subset \Lambda_p(\hat\beta_m)$ and $\det\hat{\beta}_m^\sharp(u)\neq0$ 
for all $u\in \hat{\CU}_m$ and $m\in \mathbb{N}$. Let $\{\hat{u}_m\}$ be any convergent sequence 
such that $\hat{u}_m\in \hat{\CU}_m$ for all 
$m\in\mathbb{N}$ and set $\hat{u}=\lim_{m\rightarrow\infty}\hat{u}_m$. From
$\lim_{m\rightarrow \infty}\hat{\beta}_m^\sharp(\hat{u}_m)=\hat{\beta}^\sharp(\hat{u})$ 
and since $\hat{u}_m\in\hat{\CU}_m$, it follows that $\mathrm{Index}\, \hat{\beta}^\sharp(\hat{u})<n-p$. 
Note that
$$
{\hat\beta}^\sharp(\hat u)
=\<\hat u,\hat\xi\>{\hat\beta}^\sharp(\hat\xi).
$$

We claim that $\<\hat u,\hat \xi\>=0$. 
Suppose to the contrary that $\<\hat u,\hat\xi\>>0$. If $\hat \mu>0$ (respectively, $\hat \mu<0$), 
then $\ind{\hat\beta}^\sharp(\hat u)\leq p$ (respectively, $\ind{\hat\beta}^\sharp(\hat u)\geq n-p$). Hence we have 
$\ind {\hat\beta}_m^\sharp(\hat u_m)\leq p$ (respectively, $\ind {\hat\beta}_m^\sharp(\hat u_m)\geq n-p$) 
for $m$ large enough, and that is a contradiction since 
${\hat u}_m\in \Lambda_p(\hat\beta_m)$. Similarly, we reach a contradiction if $\<\hat u,\hat\xi\><0$.
Thus we conclude that $\<\hat u,\hat \xi\>=0$.

In fact, we proved that $\langle\lim_{m\rightarrow \infty}\hat{u}_m,\hat{\xi}\rangle=0$ for any 
convergent sequence $\{\hat{u}_m\}$ such that $\hat{u}_m\in \hat{\CU}_m$ for all $m\in\N$.
We may choose convergent sequences 
$\{\hat{u}_m^{(1)}\},\dots,\{\hat{u}_m^{(k)}\}$ in $\hat{\CU}_m$ 
such that $\hat{u}_m^{(1)}, \dots, \hat{u}_m^{(k)}$ span $W$ for all $m\in \mathbb{N}$. 
We have $\langle\lim_{m\rightarrow \infty}\hat{u}_m^{(a)},\hat{\xi}\rangle=0$ for all $a\in\{1,\dots,k\}$. 
Moreover, we may write $\hat\beta_m=\sum_{a=1}^k \gamma_m^{(a)} \hat{u}_m^{(a)}$, 
with $\gamma_m^{(a)}\in \mathrm{Sym}(V\times V,\R)$ being bounded sequences. Letting $m\to \infty$ in 
the following
$$
\<\hat\beta_m(\cdot,\cdot),\hat\xi\>=\sum_{a=1}^k\gamma_m^{(a)}(\cdot,\cdot)\<\hat{u}_m^{(a)},\hat\xi\>,
$$
we obtain $\<\hat\beta(\cdot,\cdot),\hat\xi\>=0$. This yields $\hat \beta=0$, which is a contradiction. 

\smallskip
{\it Case 2}. Suppose that $\hat{\beta}$ is traceless. Given that $\hat{\beta}$ satisfies equality in 
inequality \eqref{ineq2}, Proposition \ref{lem2} yields that for all $\hat u\in \Sf^{k-1}\subset W$ 
the endomorphism ${\hat\beta}^\sharp(\hat u)$ has at most two distinct eigenvalues 
$\hat\lambda(\hat u)$ and $\hat \mu(\hat u)$ with multiplicities $p$ and $n-p$, respectively.

Since $\beta_m\in\Sigma_{n,k,p}$, it follows that there exists an open subset 
$\hat{\CU}_m\subset \Sf^{k-1}\subset W$ such that $\hat{\CU}_m\subset \Lambda_p(\hat\beta_m)$ 
and $\det\hat{\beta}_m^\sharp(u)\neq0$ for all $u\in \hat{\CU}_m$ and 
$m\in \mathbb{N}$. Let $\{\hat{u}_m\}$ be any convergent sequence such that $\hat{u}_m\in \hat{\CU}_m$ 
for all $m\in\mathbb{N}$ and set $\hat{u}=\lim_{m\rightarrow\infty}\hat{u}_m$. Since 
$\lim_{m\rightarrow \infty}\hat{\beta}_m^\sharp(\hat{u}_m)=\hat{\beta}^\sharp(\hat{u})$ 
and $\hat{u}_m\in\hat{\CU}_m$, it follows that $\mathrm{Index}\, \hat{\beta}^\sharp(\hat{u})<n-p$. 
We claim that $\hat \mu(\hat u)=0$. Observe that 
$$
p\hat\lambda(\hat u)+(n-p)\hat\mu(\hat u)=0.
$$
If $\hat\mu(\hat u)>0$ (respectively, $\hat\mu(\hat u)<0$), then $\ind {\hat\beta}^\sharp(\hat u)\leq p$ 
(respectively, $\ind {\hat\beta}^\sharp(\hat u)\geq n-p$). Thus 
$\ind {\hat\beta}_m^\sharp(\hat u_m)\leq p$ (respectively, $\ind {\hat\beta}_m^\sharp(\hat u_m)\geq n-p$) 
for $m$ large enough, and this contradicts the fact that $\hat{u}_m\in \Lambda_p(\hat\beta_m)$.
Hence $\hat \mu(\hat u)=0$ and Proposition \ref{lem2} yields ${\hat\beta}^\sharp(\hat u)=0$.

In fact, we proved that 
\be\label{limbeta}
\lim_{m\rightarrow \infty}{\hat \beta}_m^\sharp(\hat{u}_m)=0
\ee 
for any convergent sequence $\{\hat{u}_m\}$ such that $\hat{u}_m\in \hat{\CU}_m$ for all $m\in\N$.
We may choose convergent sequences $\{\hat{u}_m^{(1)}\},\dots,\{\hat{u}_m^{(k)}\}$ in $\hat{\CU}_m$ 
such that $\hat{u}_m^{(1)}, \dots, \hat{u}_m^{(k)}$ span $W$ for all $m\in \mathbb{N}$. 

Using the Gram-Schmidt process we obtain sequences $\{\xi_m^{(1)}\},\dots,\{\xi_m^{(k)}\}$ such that 
$$
\xi_m^{(i)} \in \span\{\hat u_m^{(1)},\dots,\hat u_m^{(i)}\}\;\;\text{for each}\:\: 1\leq i \leq k.
$$
Thus $\hat\xi_m^{(a)}=\sum_{\ell=1}^a x_m^{(\ell)} {\hat u}_m^{(\ell)}$, 
where $\{x_m^{(\ell)}\}$ are convergent sequences for each $1\leq \ell\leq k$,
and consequently we have
\be\label{beta1}
\hat\beta_m(\cdot,\cdot) = \sum_{a=1}^k\<{\hat\beta}_m^\sharp({\hat\xi}_m^{(a)})\cdot,\cdot\>{\hat\xi}_m^{(a)}
\ee
and ${\hat\beta}_m^\sharp({\hat\xi}_m^{(a)})= \sum_{\ell=1}^a x_m^{(\ell)}{\hat\beta}_m^\sharp({\hat u}_m^{(\ell)})$. 
Using \eqref{limbeta}, we obtain 
$$
\lim_{m\rightarrow \infty}{\hat\beta}_m^\sharp({\hat\xi}_m^{(a)})=
\lim_{m\rightarrow \infty}\sum_{\ell=1}^a x_m^{(\ell)}{\hat\beta}_m^\sharp({\hat u}_m^{(\ell)})=0 
$$ 
for all $1\leq a\leq k$. From \eqref{beta1} and the above, it follows that $\lim_{m\rightarrow \infty} \hat\beta_m=0$. 
Thus $\hat\beta=0$, which is a contradiction, and this completes the proof of Claim I. 

\medskip
Thus we may assume that $\lim_{m\rightarrow \infty}\beta_m=\beta_\infty \in\mathrm{Sym}(V\times V,W)$. 
We first argue that $\beta_\infty\neq0$. Indeed, if $\beta_\infty=0$, then $\beta_\infty^\sharp(u)=0$ 
for every $u\in \Sf^{k-1}$. Since $\beta_m\in \Sigma_{n,k,p}$, there exists $\xi_m\in \Lambda_p(\beta_m)$ 
such that 
\begin{equation}\label{CrEq6}
\left|\det\beta_m^\sharp(\xi_m)\right|\Vol(\Omega(\beta_m))=1\; \; \text{for all}\; \; m\in\mathbb{N}.
\end{equation}
We may assume that the sequence $\{\xi_m\}$ converges to some $\xi\in\Sf^{k-1}$. Hence 
we have $\lim_{m\rightarrow \infty}\beta_m^\sharp(\xi_m)=\beta_\infty^\sharp(\xi)=0$, which 
contradicts (\ref{CrEq6}). 

\medskip
{\it Claim II: The bilinear form $\beta_\infty$ is not a zero of $\phi_p$.} Arguing by contradiction, 
suppose that $\phi_p(\beta_\infty)=0$, that is $\beta$ satisfies equality in inequality 
\eqref{ineq2}. To reach a contradiction and according to Proposition \ref{lem2}, 
we distinguish two cases.

\smallskip
{\it Case 3}. Suppose that $\mathsf H_{\beta_\infty}\neq 0$. Hence 
$\beta_\infty=\< {\beta}_\infty^\sharp(\xi) \cdot , \cdot \>\, \xi$, where 
$\xi=\mathsf H_{\beta_\infty}/\|\mathsf H_{\beta_\infty}\|$ and ${\beta}_\infty^\sharp(\xi)$ has 
at most two distinct eigenvalues $\lambda$ and $\mu$ with multiplicities $p$ and $n-p$, 
respectively. Clearly $\mu\neq0$, since otherwise $\beta_\infty=0$, which is a contradiction.

Given that $\beta_m\in \Sigma_{n,k,p}$, there exists an open subset $\CU_m\subset \Sf^{k-1}\subset W$ 
such that $\CU_m\subset \Lambda_p(\beta_m)$ and $\det{\beta}_m^\sharp(u)\neq0$ for each 
$u\in \CU_m$ and $m\in \mathbb{N}$. Let $\{u_m\}$ be any convergent sequence such that 
$u_m\in \CU_m$ for all $m\in\mathbb{N}$ and set $u=\lim_{m\rightarrow\infty}u_m$. Since 
$\lim_{m\rightarrow \infty}{\beta}_m^\sharp(u_m)={\beta}_\infty^\sharp(u)$ 
and $u_m\in\CU_m$, it follows that $\mathrm{Index}\, {\beta}_\infty^\sharp(u)<n-p$. 
Observe that
$$
{\beta}_\infty^\sharp(u)
=\<u,\xi\>{\beta}_\infty^\sharp(\xi).
$$

We claim that $\< u, \xi\>=0$. Suppose to the contrary that $\<u,\xi\>>0$. If $\mu>0$ 
(respectively, $\mu<0$), then $\ind{\beta}_\infty^\sharp(u)\leq p$ 
(respectively, $\ind{\beta}_\infty^\sharp(u)\geq n-p$) and consequently 
$\ind{\beta}_m^\sharp(u_m)\leq p$ (respectively, $\ind{\beta}_m^\sharp(u_m)\geq n-p$) for $m$ 
large enough, and that is a contradiction since $u_m\in \Lambda_p(\beta_m)$.
A contradiction is reached in a similar manner if $\<u,\xi\><0$. 
Thus we conclude that $\<u,\xi\>=0$. 

In fact, we proved that $\langle\lim_{m\rightarrow \infty}u_m,\xi\rangle=0$ for any convergent 
sequence $\{u_m\}$ such that $u_m\in \CU_m$ for all $m\in\N$. We may choose convergent 
sequences $\{u_m^{(1)}\},\dots,\{u_m^{(k)}\}$ in $\CU_m$ 
such that $u_m^{(1)}, \dots, u_m^{(k)}$ span $W$ for all $m\in \mathbb{N}$. 
We have $\langle\lim_{m\rightarrow \infty}u_m^{(a)},\xi\rangle=0$ for all $a\in\{1,\dots,k\}$. 
Moreover, we may write $\beta_m=\sum_{a=1}^k \gamma_m^{(a)} u_m^{(a)}$, 
with $\gamma_m^{(a)}\in \mathrm{Sym}(V\times V,\R)$ being bounded. Letting $m\to\infty$ 
in the equality
$$
\<\beta_m(\cdot,\cdot),\xi\>=\sum_{a=1}^k\gamma_m^{(a)}(\cdot,\cdot)\<u_m^{(a)},\xi\>,
$$
it follows that $\<\beta_\infty(\cdot,\cdot),\xi\>=0$. 
Thus $\beta_\infty=0$, which is a contradiction. 

\smallskip
{\it Case 4}. Suppose that $\beta_\infty$ is trace free. Given that $\beta_\infty$ satisfies equality in 
inequality \eqref{ineq2}, Proposition \ref{lem2} yields that for all $u\in \Sf^{k-1}\subset W$ 
the endomorphism ${\beta}_\infty^\sharp(u)$ has at most two distinct eigenvalues $\lambda(u)$ 
and $\mu(u)$ with multiplicities $p$ and $n-p$, respectively.

Since $\beta_m\in \Sigma_{n,k,p}$, there exists an open subset $\CU_m\subset \Sf^{k-1}\subset W$ 
such that $\CU_m\subset \Lambda_p(\beta_m)$ and $\det{\beta}_m^\sharp(u)\neq0$ for all 
$u\in \CU_m$ and $m\in \mathbb{N}$. Let $\{u_m\}$ be any convergent sequence such that 
$u_m\in\CU_m$ for all $m\in\mathbb{N}$ with $u=\lim_{m\rightarrow\infty}u_m$. Since 
$\lim_{m\rightarrow \infty}{\beta}_m^\sharp(u_m)={\beta}_\infty^\sharp(u)$ 
and $u_m\in\CU_m$, it follows that $\mathrm{Index}\, {\beta}_\infty^\sharp(u)<n-p$. 

We claim that $\mu(u)=0$. Observe that 
$$
p\lambda(u)+(n-p)\mu(u)=0.
$$
If $\mu(u)>0$ (respectively, $\mu(u)<0$), then $\ind {\beta}_\infty^\sharp(u)\leq p$ 
(respectively, $\ind {\beta}_\infty^\sharp(u)\geq n-p$).
Hence $\ind {\beta}_m^\sharp(u_m)\leq p$ (respectively, $\ind {\beta}_m^\sharp(u_m)\geq n-p$) 
for $m$ large enough, and that is a contradiction since $u_m\in \Lambda_p({\beta}_m)$. 
Thus we have $\mu(u)=0$, which implies that ${\beta}_\infty^\sharp(u)=0$.

In fact, we proved that 
\be\label{limbeta.2}
\lim_{m\rightarrow \infty}{\beta}_m^\sharp(u_m)=0
\ee 
for every convergent sequence $\{u_m\}$ such that $u_m\in \CU_m$ for all $m\in\N$.
We may choose convergent sequences 
$\{u_m^{(1)}\},\dots,\{u_m^{(k)}\}$ in $\CU_m$ 
such that $u_m^{(1)}, \dots, u_m^{(k)}$ span $W$ for all $m\in \mathbb{N}$. 

Using the Gram-Schmidt process, we obtain sequences $\{\xi_m^{(1)}\},\dots,\{\xi_m^{(k)}\}$ such that 
$$
\xi_m^{(i)} \in \span\{u_m^{(1)},\dots,u_m^{(i)}\}\;\;\text{for each}\;\; 1\leq i\leq k.
$$
Thus $\xi_m^{(a)}=\sum_{\ell=1}^a x_m^{(\ell)} u_m^{(\ell)}$, 
where $\{x_m^{(\ell)}\}$ are convergent sequences for every $1\leq \ell\leq k$.
Clearly we have
\be\label{beta1.2}
\beta_m(\cdot,\cdot) = \sum_{a=1}^k\<{\beta}_m^\sharp({\xi}_m^{(a)})\cdot,\cdot\>{\xi}_m^{(a)}
\ee
and therefore ${\beta}_m^\sharp({\xi}_m^{(a)})= \sum_{\ell=1}^a x_m^{(\ell)}{\beta}_m^\sharp(u_m^{(\ell)})$. 
Using \eqref{limbeta.2}, it follows that
$$
\lim_{m\rightarrow \infty}{\beta}_m^\sharp({\xi}_m^{(a)})=
\lim_{m\rightarrow \infty}\sum_{\ell=1}^a x_m^{(\ell)}{\beta}_m^\sharp(u_m^{(\ell)})=0
$$
for all $1\leq a\leq k$. From \eqref{beta1.2} and the above it follows that $\beta_\infty=0$, 
which is a contradiction. This completes the proof of Claim II. 

\smallskip
Thus the function $\phi_p$ attains a positive minimum 
$\varepsilon(n,k,p)=\phi_p(\beta_\infty)$ on $\Sigma_{n,k,p}$, which depends only on $n,k$ and $p$. 
Let $\beta$ be an arbitrary bilinear form in $\mathrm{Sym}(V\times V, W)$. Suppose that 
$\psi_p(\beta)\neq 0$ and set $\gamma=\beta/(\psi_p(\beta))^{1/n}$. Since $\gamma\in\Sigma_{n,k,p}$, 
we have $\phi_p(\gamma)\geq \varepsilon(n,k,p)$ and the desired inequality follows from the 
homogeneity of the function $\phi_p$.
\qed

\section{Submanifolds satisfying the pinching condition}

Let $M^n$ be an oriented Riemannian manifold of dimension $n$. For each integer $0\leq p\leq n$, 
the Hodge-Laplace operator acting on differential $p$-forms is defined by 
$$
\Delta=d\delta+\delta d: \Omega^p(M^n)\to \Omega^p(M^n),
$$
where $d$ and $\delta$ are the differential and the co-differential operators, respectively. 
When $p=0$, the Hodge-Laplace operator is just the Laplace-Beltrami operator acting on 
$0$-forms, i.e., scalar functions.

We recall the following well known result (cf. \cite[Prop. 3]{Savo14}).

\begin{proposition}\label{bochner}
Let $M^n$ be a compact, oriented Riemannian manifold and let $1\leq p\leq n/2$.
The following assertions hold:
\begin{enumerate}[(i)]
\item If $\B^{[p]}\geq 0$, then $b_p(M^n)\leq \binom n p$.
\item If $\B^{[p]}\geq 0$ and the strict inequality holds at some point, then 
$H^p(M^n;\R)=H^{n-p}(M^n;\R)=0$.
\item If $\B^{[p]}\geq 0$ and $H^p(M^n;\R)\neq0$, then every harmonic $p$-form is parallel. 
In particular, $M^n$ supports a nontrivial parallel $p$-form.
\end{enumerate}
\end{proposition}

\proof
Let $\omega\in \Omega^p(M^n)$. From the Bochner-Weitzenb\"ock formula \eqref{boch.form}, we obtain 
\be\label{boch.formula}
\<\Delta \omega,\omega\>=\|\n\omega\|^2+\<\B^{[p]}\omega,\omega \>+\frac{1}{2}\, \Delta \|\omega\|^2.
\ee
If $\omega$ belongs to the space $\mathscr H^p(M^n)$ of harmonic $p$-forms, then using 
the assumption, we have 
$$
0\leq\int_{M^n} \<\B^{[p]}\omega,\omega \>dM=-\int_{M^n}\|\n\omega\|^2dM. 
$$
Thus $\omega$ is parallel and therefore
\begin{eqnarray}
b_p(M^n) \!\!\!&=&\!\!\! \dim H^p(M^n; \R)=\dim \mathscr H^p(M^n) \nonumber \\
\!\!\!&\leq&\!\!\! \dim \Lambda^p(T_x^*M^n)=\binom n p. \nonumber 
\end{eqnarray}
This proves (i). For the proofs of (ii) and (iii) see \cite[Prop. 3]{Savo14}.
\qed

\medskip

Let $f\colon M^n\to \tilde M^{n+k}$ be an isometric immersion whose second fundamental form 
$\a_f$ is viewed as a section of the vector bundle $\mathrm{Hom}(TM\times TM,N_f M)$, where 
$N_f M$ is the normal bundle. For each unit normal vector field $\xi\in \Gamma(N_fM)$, the associated 
shape operator $A_{\xi}$ is given by 
$$
\<A_{\xi} X,Y\>=\langle \alpha_f(X,Y),\xi\rangle,\: \: X,Y\in TM.
$$
 
Savo \cite[Th. 1]{Savo14} proved that the Bochner operator $\B^{[p]}$ splits as
\be\label{bosplit}
\B^{[p]}=\B^{[p]}_{\rm{res}}+\B^{[p]}_{\rm{ext}},
\ee
where $\B^{[p]}_{\rm{res}}$ depends on the geometry of the ambient manifold $\tilde M^{n+k}$
restricted to the submanifold $f(M^n)$ and satisfies 
\be\label{bo1}
\B^{[p]}_{\rm{res}}\geq p(n-p)c, 
\ee
where $c$ is a lower bound for the curvature operator of $\tilde M^{n+k}$. In particular, if 
$\tilde M^{n+k}$ has constant sectional curvature $c$, then 
\be\label{bo2}
\B^{[p]}_{\rm{res}}=p(n-p)c.
\ee
The second operator $\B^{[p]}_{\rm{ext}}$ is given 
explicitly in terms of the second fundamental form of the immersion $f$. 
For each unit vector field $\xi\in \Gamma(N_fM)$, we define 
the endomorphism 
$$
T_{A_\xi}^{[p]}= (\tr A_\xi)A_\xi^{[p]}-A_\xi^{[p]}\circ A_\xi^{[p]},
$$
where 
$$
A_\xi^{[p]}\colon \Omega^p(M^n)\to \Omega^p(M^n)
$$
is the self-adjoint extension of the shape operator $A_\xi$ in the direction $\xi$ defined by 
$$
A_\xi^{[p]}\omega(X_1,\dots,X_p)=\sum_{i=1}^p \omega(X_1,\dots,A_\xi X_i,\dots,X_p),\; X_1,\dots,X_p\in TM.
$$
Then $\B^{[p]}_{\rm{ext}}$ is given by (see \cite{Savo14})
$$
\B_{\rm{ext}}^{[p]}=\sum_{i=1}^k T_{A_{\xi_i}}^{[p]},
$$
where $\xi_1,\dots, \xi_k$ is an orthonormal frame of the normal bundle of $f$.

\medskip
The next proposition provides a sharp estimate for the Bochner operator in terms of the 
second fundamental form in arbitrary codimension and gives its structure at points where 
equality holds. Inequality \eqref{ineq4} was proved by Savo \cite[Prop. 15]{Savo14} for 
hypersurfaces. Note that, since $\B^{[1]}$ is nothing but the Ricci tensor, inequality \eqref{ineq4} for 
$p=1$ reduces to the inequality due to Leung \cite{leung92}. We recall that the traceless part of the 
second fundamental form is given by $\Phi=\a_f-\<\cdot,\cdot\>\mathcal H$. 

\begin{proposition}\label{lem3}
Let $f\colon M^n\to \tilde M^m, n\geq 3$, be an isometric immersion. If the 
curvature operator of $\tilde M^m$ is bounded from below by a constant $c$, then the Bochner 
operator of $M^n$ satisfies pointwise the inequality
\be\label{ineq4}
\mathop{\min_{\omega\in \Omega^p(M^n)}}_{\|\omega\|=1}\<\B^{[p]}\omega,\omega\>
\geq \frac{p(n-p)}{n} \big(n(H^2+c) -\frac{n(n-2p)}{\sqrt{np(n-p)}}\, H\|\Phi\|-\|\Phi\|^2\big).
\ee
If equality holds in \eqref{ineq4} at a point $x\in M^n$, then the following hold:
\begin{enumerate}[(i)]
\item The shape operator $A_\xi(x)$ has at most two distinct 
eigenvalues with multiplicities $p$ and $n-p$ for every unit vector $\xi\in N_fM(x)$. If in addition 
$p<n/2$ and the eigenvalue of multiplicity $n-p$ vanishes, then $A_\xi(x)=0$.
\item If $H(x)\neq 0$ and $p<n/2$, then $\mathrm{Im}\, \a(x)=\span\left\{\mathcal H(x)\right\}$.
\end{enumerate}
\end{proposition}

\proof
It is a direct consequence of Proposition \ref{lem2}, \eqref{bosplit} and \eqref{bo1}.
\qed

\medskip

{\noindent{\it Proof of Theorem \ref{mainthm1}.}}
We observe that
\be\label{ar}
a(n,p,H,c)=\left(r(n,p,H,c)\right)^2+nH^2,
\ee
where $r(n,p,H,c)$ is the largest root of the quadratic polynomial
\be\label{P}
P(t;p,H,c)= t^2+\frac{n(n-2p)}{\sqrt{np(n-p)}}\, H t-n(H^2+c).
\ee
Using that $\|\Phi\|^2=S-nH^2$, our pinching assumption $S\leq a(n,p,H,c)$ turns out to be equivalent to $\|\Phi\|\leq r(n,p,H,c)$. 
Since
$$
0\leq r(n,1,H,c)\leq \cdots\leq r(n,\left[n/2\right],H,c),
$$
we have $\|\Phi\|\leq r(n,i,H,c),$ 
or equivalently 
$$
P(\|\Phi\|;i,H,c)\leq 0 \;\;\text{for each}\;\; p\leq i\leq n/2.
$$ 
Proposition \ref{lem3} implies that 
$\B^{[i]}\geq0$ for all $p\leq i\leq n/2$. Now, part (i) of the theorem follows from Proposition 
\ref{bochner}(i) and the Poincar\'e duality, while part (ii) follows from Propositions \ref{bochner}(ii) 
and \ref{lem3}. Part (iii) follows from part (ii) and Proposition \ref{bochner}(iii). \qed

\medskip

{\noindent{\it Proof of Theorem \ref{thm2}.}} Let $\omega \in\mathscr H^p(M^n)$ be a harmonic
 form of constant length. Equation \eqref{boch.formula} yields
$$
\<\B^{[p]}\omega,\omega\>=-\|\n\omega\|^2\leq 0.
$$
From Proposition \ref{lem3}, it follows that 
$$
n(H^2+c) -\frac{n(n-2p)}{\sqrt{np(n-p)}}\, H\|\Phi\|-\|\Phi\|^2\leq 0.
$$
Hence we have $\|\Phi\| \geq r(n,p,H,c)$, where $r(n,p,H,c)$ is the largest root of the polynomial \eqref{P}. 
Since $\|\Phi\|^2=S-nH^2$,  using \eqref{ar} we obtain $S\geq a(n,p,H,c)$.
\qed

\medskip
We recall that an isometric immersion $f\colon M^n\to \Q_c^{n+k}$, into a space form of constant sectional 
curvature $c$, is said to reduce codimension to $m <k$ if there exists a totally geodesic submanifold 
$\Q_c^{n+m}$ of $\Q_c^{n+k}$ such that $f(M)\subset \Q_c^{n+m}$. We need the following proposition 
(a related result was given in \cite{saco12} for hypersurfaces in spheres).

\begin{proposition}\label{cod}
Let $f\colon M^n\to \Q_c^{n+k}, n\geq 3$, be an isometric immersion into a nonflat space form that 
satisfies \eqref{1} for an integer $1\leq p< n/2$ and $H^2+c\geq0$. Suppose 
that $M^n$ supports a nontrivial harmonic $p$-form (in particular, a parallel $p$-form) with constant 
length. If $H>0$ on $M^n$, then the codimension of $f$ is reduced to one. Moreover, the hypersurface 
$f\colon M^n\to \Q_c^{n+1}$ has at most two distinct principal curvatures $\lambda$ and $\mu$ 
of multiplicities $p$ and $n-p$, respectively, such that
$\lambda\mu+c=0.$
\end{proposition}

\proof
Let $\omega\in{\mathscr H}^p(M^n)$ be a nontrivial harmonic $p$-form of constant length. 
Equation \eqref{boch.formula} yields $\<\B^{[p]}\omega,\omega \>\leq0$. Theorem 
\ref{thm2} 
implies that $S=a(n,p,H,c)$ on $M^n$, which is equivalent to
$$
n(H^2+c) -\frac{n(n-2p)}{\sqrt{np(n-p)}}\, H\|\Phi\|-\|\Phi\|^2=0.
$$
Thus equality holds in \eqref{ineq4} at each point. According to Proposition \ref{lem3}(ii), the second 
fundamental form of $f$ satisfies
$$
\alpha_f(X,Y)=\<A_\xi X,Y\>\xi \;\;\text{for all}\;\; X,Y\in TM,
$$
where $\xi=\mathcal H/H$ and the shape operator $A_\xi$ has at most two distinct eigenvalues 
$\lambda$ and $\mu$ of multiplicities $p$ and $n-p$, respectively. The Gauss equation is 
written as
$$
R(X,Y)Z=c\left(\<Y,Z\>X-\<X,Z\>Y\right)+\<A_\xi Y,Z\>A_\xi X-\<A_\xi X,Z\>A_\xi Y,
$$
where the curvature tensor $R$ is defined by 
$$
R(X,Y)=[\n_X,\n_Y]-\n_{[X,Y]}, \;X,Y \in \mathcal X(M^n).
$$
Then we find that 
\be\label{Romega}
R(X,Y)\omega= cX^\#\wedge i_Y\omega-cY^\#\wedge i_X \omega
+(A_\xi X)^\#\wedge i_{A_\xi Y}\omega-(A_\xi Y)^\#\wedge i_{A_\xi X}\omega.
\ee
Here $i_X \omega$ is the interior multiplication of the form $\omega$ by the vector field $X$ and 
$X^\#$ is the 1-form dual to $X$. Let $X_1,\dots,X_n$ be an arbitrary local orthonormal frame with 
dual frame $\theta_1,\dots,\theta_n$ and set
$$
\Psi_{ij}=\theta_i\wedge i_{X_j} \omega-\theta_j\wedge i_{X_i}\omega, \;\;1\leq i\leq p<j\leq n.
$$
We claim that $\Psi_{ij}\neq 0$. In fact, since $\omega$ is parallel and nontrivial, we may assume 
after reordering that 
$\omega(X_1\dots,X_p)\neq0$. Thus we have
$$
\Psi_{ij}(X_1,\dots,\hat X_i ,\dots,X_p, X
_j)=\pm\omega(X_1\dots,X_p)\neq0.
$$

Now we choose the local orthonormal frame $X_1,\dots,X_n$ such that 
$$
A_\xi X_i=\lambda X_i, \; 1\leq i\leq p
\; \; \text{and}\; \;
 A_\xi X_j=\mu X_j,\; p+1\leq j\leq n.
$$
Since $\omega$ is parallel, we have $R(X_i,X_j)\omega=0$ and equation \eqref{Romega} yields 
$(\lambda\mu+c)\Psi_{ij}=0$, or equivalently $\lambda\mu+c=0$.

For each $\eta\in \Gamma(N_fM)$ perpendicular to $\xi$, 
the Codazzi equation yields 
$$
A_{\n_{X_i}^\perp \eta} X_j=A_{\n_{X_j}^\perp \eta} X_i,\; 1\leq i\leq p<j\leq n,
$$ 
or equivalently 
$$
\mu\<{\n_{X_i}^\perp \eta},\xi\> X_j=\lambda\<{\n_{X_j}^\perp \eta},\xi\> X_i,\; 1\leq i\leq p<j\leq n.
$$ 
This implies that
$$
\<{\n_{X_i}^\perp \eta},\xi\>=\<{\n_{X_j}^\perp \eta},\xi\>=0
$$ 
for all $1\leq i\leq p$ and $p+1\leq j\leq n$. Therefore the first normal bundle
$N^1_fM=\mathrm{Im}\, \a_f$ is a parallel subbundle of the normal bundle of constant rank one. 
Thus the codimension of $f$ is reduced to one (see for instance \cite[Prop. 2.1]{DT}). \qed

\medskip

{\noindent{\it Proof of Theorem \ref{thm.b}.}} 
Suppose that $M^n$ is not a real homology sphere and let $p$ be an integer such that $b_p(M^n)>0$ 
with $1\leq p\leq n/2$. Since $S\leq n=a(n,p,0,1)$, it follows from Theorem \ref{mainthm1}(iii) 
that $S=n$ on $M^n$ and that every harmonic form $\omega\in {\mathscr H}^p(M^n)$ is parallel. 
Let $\omega$ be a non vanishing parallel form. Equation \eqref{boch.formula} yields 
$\<\B^{[p]}\omega,\omega \>=0$. Hence equality holds in inequality \eqref{ineq4} at every point. 
Proposition \ref{lem3}(i) implies that the second fundamental form 
in every normal direction has at most two distinct principal curvatures of multiplicities $p$ and $n-p$.

If $m=n+1$, then it follows from \cite{CDCK70} that 
$f(M^n)$ is isometric to the Clifford torus $\mathbb T^n_p(\sqrt{p/n})$. 
If $m>n+1$, then it follows from \cite[Th. 1.5]{jome84} that $f$ is a standard embedding of a 
projective plane over the complex $\C$, quaternions $\mathbb H$, or Cayley numbers $\mathbb O$. 
In other words, $f$ is one of the standard minimal isometric embeddings
$$
\psi_1 \colon \CP^2\to \Sf^7,\; \psi_2 \colon \HP^2\to \Sf^{13},\; \psi_3\colon {\mathbb{OP}}^2\to \Sf^{25}.
$$
The holomorphic curvature of $\CP^2$ is $4/3$, the $\mathbb H$-sectional curvature of $\HP^2$ is 
$4/3$ and the ${\mathbb O}$-sectional curvature of the Cayley plane is $4/3$
 (see \cite[p. 780]{gauchman84} or \cite{ish74}). It can be easily checked that 
 $\psi_1 \colon \CP^2\to \Sf^7$ is the only one of the above satisfying $S=n$.

Now suppose that $M^n$ is a real homology sphere. It remains to prove that $M^n$ admits a metric of positive Ricci curvature.
Using our assumptions, it follows from inequality \eqref{ineq4}  for $p=1$ that the Ricci curvature of $M^n$ is nonnegative
everywhere. 

We claim that there is a point where the Ricci curvature is positive for every tangent direction.
Suppose to the contrary that at each point there exists a tangent direction where the Ricci curvature vanishes.
This implies that equality holds in inequality \eqref{ineq4} for $p=1$ at every point. Proposition \ref{lem3}(i) implies that the 
second fundamental form in every normal direction has at most two distinct principal curvatures, one with multiplicity 
$n-1$. Arguing as before, it follows from \cite{CDCK70} and \cite[Th. 1.5]{jome84} that either $f(M^n)$ 
is isometric to the Clifford torus $\mathbb T^n_1(\sqrt{1/n})$, or $f$ is the standard embedding 
$\psi\colon \CP^2_{4/3}\to \Sf^7$ of the complex projective plane of constant holomorphic curvature $4/3$.
Hence either $M^n$ is diffeomorphic to $\Sf^1\times\Sf^{n-1}$ (cf. \cite[Th. 1]{fz}), or diffeomorphic to $\CP^2$.
This clearly contradicts the fact that $M^n$ is a real homology sphere. 

Thus there exists a point where the Ricci curvature is positive for every tangent direction.
It follows from Aubin \cite{A} that $M^n$ carries a metric with positive Ricci curvature. 
The Bonnet-Myers theorem implies that the fundamental group of $M^n$ is finite. 
\qed

\medskip

{\noindent{\it Proof of Corollary \ref{cor}.}}
Suppose that $n=3$ and $f(M^3)$ is not isometric to the Clifford torus $\mathbb T^3_1(\sqrt{1/3})$. 
Theorem \ref{thm.b} implies that $M^3$ carries a metric with positive Ricci curvature. It follows from 
\cite[Th. 1.1]{ham} that $M^3$ is diffeomorphic to a spherical space form.

Now assume that $n=4$ and the minimal submanifold is neither a Clifford torus, nor the standard embedding 
$\psi$ of the complex projective plane. Theorem \ref{thm.b} implies that $M^4$ is a real homology 
sphere and its universal cover $\tilde M^4$ is compact. 
We claim that $\tilde M^4$ is a real homology sphere. If otherwise, 
we consider the minimal isometric immersion given by $\tilde f= f \circ \pi$, where 
$\pi\colon\tilde M^4 \to M^4$ is the covering map. It is obvious that the squared length 
of its second fundamental form satisfies $\tilde S\leq n$. Then Theorem \ref{thm.b} would imply that either 
$f(M^4)$ has to be isometric to a minimal Clifford torus, or $\tilde f$ is the standard isometric embedding 
of the complex projective plane. This is clearly a contradiction.

From the Poincar\'e duality and the universal coefficient theorem (see \cite[Cor. 4, p. 244]{sp}) 
it follows that the torsion subgroups of $H_i(\tilde M^4;\Z)$ and $H_{3-i}(\tilde M^4;\Z)$ are 
isomorphic for $i=1,2,3$. Since $\tilde M^4$ is both simply connected and a real homology sphere, 
we have that it is also a homology sphere over the integers. Therefore $\tilde M^4$ is a homotopy 
sphere. By the proof of the Poincar\'e conjecture for $n=4$ due to Freedman \cite{Fr}, the universal cover $\tilde M^4$ is 
homeomorphic to $\Sf^n$.
\qed

\medskip

{\noindent{\it Proof of Theorem \ref{thm.a}.}} Suppose that $b_p(M^n)>0$ and let 
$\omega\in{\mathscr H}^p(M^n)$ be a non vanishing harmonic form. Theorem \ref{mainthm1} (iii) 
implies that the form $\omega$ is parallel and $S=a(n,p,H,1)$
on $M^n$, or equivalently 
$$
n(H^2+1) -\frac{n(n-2p)}{\sqrt{np(n-p)}}\, H\|\Phi\|-\|\Phi\|^2=0.
$$
Then Proposition \ref{lem3} gives $\<\B^{[p]}\omega,\omega \>\geq0$. On the other hand, since 
$\omega$ is parallel, it follows from equation \eqref{boch.formula} that $\<\B^{[p]}\omega,\omega \>=0$.
Hence equality holds in inequality \eqref{ineq4} at every point, and Proposition \ref{lem3}(i) implies 
that the second fundamental form of $f$ has at most two distinct principal curvatures of multiplicities 
$p$ and $n-p$ in each normal direction. 

Let $v\in\R^{m+1}$ be a
vector such that the height function $h\colon M^n\to \R$ defined by $h=\<g,v\>$ is a Morse function, 
where $g$ is the isometric immersion $g=j\circ f$, and $j\colon \Sf^m\to \R^{m+1}$ is the inclusion. 
The Hessian of $h$ is given by
$$
\mathrm{Hess}\, h(X,Y)=\< \a_g(X,Y),v\>,\;X,Y\in TM.
$$
Obviously, the second fundamental form of $g$ has at most two distinct principal curvatures of 
multiplicities $p$ and $n-p$ in every normal direction. Hence the index of each nondegenerate 
critical point of $h$ is $0,p,n-p$ or $n$. That the manifold $M^n$ has the homotopy type of a 
CW-complex with cells only in dimensions $0,p,n-p$ or $n$ follows from standard Morse theory 
(cf. \cite[Th. 3.5]{Milnor63} or \cite[Th. 4.10]{CE75}), and this proves part (i) of the theorem.

(ii) Our pinching assumption $S\leq a(n,p,H,1)$ turns out to be equivalent to 
$$
\|\Phi\|\leq r(n,p,H,1)\leq r(n,q,H,1).
$$
Since $b_q(M^n)>0$, Theorem \ref{mainthm1} (iii) 
implies that equality holds everywhere in the above inequalities, and consequently
$$
n(H^2+1) -\frac{n(n-2p)}{\sqrt{np(n-p)}}\, H\|\Phi\|-\|\Phi\|^2=0
$$
and 
$$
n(H^2+1) -\frac{n(n-2q)}{\sqrt{nq(n-q)}}\, H\|\Phi\|-\|\Phi\|^2=0.
$$
It follows directly from the above that the submanifold is minimal and $S=n$. Hence 
Theorem \ref{thm.b} implies that either $f(M^n)$ is isometric to the Clifford torus 
$\mathbb T^n_q(\sqrt{q/n})$, or $f$ is the standard embedding of the complex projective plane.

Hereafter, we suppose that $H>0$ on $M^n$.

(iii) Since $b_p(M^n)>0$, Theorem \ref{mainthm1} implies that there exist a nontrivial parallel 
$p$-form. Using our assumption on the mean curvature, it follows from Proposition \ref{cod} 
that the codimension of $f$ reduces to one, and the proof of part (iii) follows from \cite[Th. 9]{Savo14}.

(iv) Under our assumption, we obtain 
$$
S=a(n,p,H,1)<a(n,p+1,H,1)
$$
on $M^n$. From \cite[Th. 1]{Vlachos07}, it follows that the homology groups satisfy 
$H_i(M^n;\Z)=0$ for every $p< i< n-p$. Now suppose that $f(M^n)$ is not isometric to a torus 
as in the statement of part (iii) of the theorem. Then we have $b_p(M^n)=0$. From the 
universal coefficient theorem (see \cite[Cor. 4, p. 244]{sp}) and using $H_{n-p-1}(M^n;\Z)=0$, 
we conclude that $H^{n-p}(M^n;\Z)$ has no torsion and neither does $H_p(M^n;\Z)$ 
by Poincar\'e duality. Since $b_p(M^n)=0$, we obtain $H_p(M^n;\Z)=0$. 
\qed

\medskip

{\noindent{\it Proof of Theorem \ref{cor0}.}}
At first we observe that $M^n$ has nonnegative Ricci curvature. Indeed, our assumption $S\leq a(n,1,H,1)$ 
is equivalent to $\|\Phi\|\leq r(n,1,H,1)$, where $r(n,1,H,1)$ is the largest root of the polynomial 
\eqref{P}. Hence the right hand side of inequality \eqref{ineq4} for $p=1$ is nonnegative, and consequently 
$M^n$ has nonnegative Ricci curvature. 

If $b_1(M^n)>0$, then it follows from Theorem \ref{thm.a}(i) that $M^n$ has the homotopy type of a 
CW-complex with cells only in dimensions $0,1,n-1,n$. If $b_q(M^n)>0$ for an integer $1<q\leq n/2$, then 
Theorem \ref{thm.a}(ii) yields that $f(M^n)$ is isometric to the Clifford torus $\mathbb T^n_q(\sqrt{q/n})$ 
for an integer $1<q\leq n/2$, or $f$ is the standard embedding as in Theorem \ref{thm.b}(iii). 

Suppose now that none of the above holds. Hence, according to Theorem \ref{thm.a}, $M^n$ is a real 
homology sphere. We claim that there is a point where all Ricci curvatures are positive. If otherwise, 
equality holds in \eqref{ineq4} for $p=1$ everywhere, and consequently the second fundamental form 
has at most two distinct principal curvatures of multiplicities $1$ and $n-1$ in every normal direction. 
Then arguing as in the proof of Theorem \ref{thm.a} using Morse theory, we conclude that $M^n$ has 
the homotopy type of a CW-complex with cells only in dimensions $0,1,n-1,n$. This is a contradiction 
since this case has been excluded.

Hence there is a point where all Ricci curvatures are positive. It follows from Aubin \cite{A} that $M^n$
carries a metric of positive Ricci curvature. By Bonnet-Myers theorem, the fundamental group of $M^n$
is finite and its universal covering $\tilde M^n$ is compact.  We claim that $\tilde M^n$ is a real homology sphere.
In fact, the isometric immersion $\tilde f=f \circ \pi$ clearly satisfies our pinching condition for $p=1$, where 
$\pi\colon\tilde M^n \to M^n$ is the covering map. At first, we observe that $b_1(\tilde M^n)=0$. If otherwise, 
Theorem \ref{thm.a}(i) would imply that $\tilde M^n$ has the homotopy type of a CW-complex with cells only in 
dimensions $0,1,n-1,n$, and this is clearly a contradiction. If $b_q(\tilde M^n)>0$ for an integer $1<q\leq n/2$, 
then Theorem \ref{thm.a}(ii) yields that $f(M^n)$ is isometric to the Clifford torus $\mathbb T^n_q(\sqrt{q/n})$, 
or $f$ is the standard embedding as in Theorem \ref{thm.b}(iii). This is again a contradiction.
Thus $\tilde M^n$ is a real homology sphere. The assertion for $n=3,4$ follows 
as in the proof of Corollary \ref{cor}. This completes the proof of parts (i)-(iii).

Now suppose that $H>0$ on $M^n$. Assume that there exists a point where the Ricci curvature is positive, then 
according to Aubin \cite{A}, the manifold $M^n$ admits a metric of positive Ricci curvature. The 
Bonnet-Myers theorem implies that the fundamental group of $M^n$ is finite. Since we have
$$
S\leq a(n,1,H,1)<a(n,2,H,1)
$$
on $M^n$,  it follows from \cite[Cor. 1]{Vlachos07} that $M^n$ is homeomorphic to $\Sf^n$.

Now suppose that $M^n$ is not homeomorphic to $\Sf^n$. Thus at every point there exists a 
tangent direction where the Ricci curvature vanishes. This implies that the left hand side of inequality 
\eqref{ineq4} is nonpositive, being the right hand side nonnegative.
Hence equality holds in \eqref{ineq4} at every point. Proposition \ref{lem3}(ii) implies that 
the second fundamental form of $f$ satisfies
$$
\alpha_f(X,Y)=\<A_\xi X,Y\>\xi \;\;\text{for all}\;\; X,Y\in TM,
$$
where $\xi=\mathcal H/H$ and the shape operator $A_\xi$ has at most two distinct eigenvalues 
$\lambda$ and $\mu$ of multiplicities $1$ and $n-1$, respectively. We choose a local orthonormal 
frame $X_1,\dots,X_n$ such that 
$$
A_\xi X_1=\lambda X_1 \; \; \text{and}\; \; A_\xi X_i=\mu X_i,\; 2\leq i\leq n.
$$
For each $\eta\in \Gamma(N_fM)$ perpendicular to $\xi$, the Codazzi equation yields 
$$
\mu\<{\n_{X_1}^\perp \eta},\xi\>X_i=\lambda\<{\n_{X_i}^\perp \eta},\xi\> X_1,\; 2\leq i\leq n
$$ 
and
$$
\mu\<{\n_{X_i}^\perp \eta},\xi\> X_j=\mu\<{\n_{X_j}^\perp \eta},\xi\> X_i,\; 2\leq i\neq j \leq n.
$$ 
We claim that $\mu$ cannot vanish. Indeed, if $\mu$ vanishes at a point $x\in M^n$, then 
Proposition \ref{lem3}(i) gives that $x$ is a totally geodesic point and this contradicts our assumption on 
the mean curvature. Hence we have $\<{\n_{X_i}^\perp \eta},\xi\>=0$ for all $1\leq i\leq n$, and 
therefore the first normal bundle is a parallel subbundle of the normal bundle of constant rank one. 
Thus $f$ reduces codimension to one and can be regarded as a hypersurface in a totally 
geodesic $\Sf^{n+1}$ of $\Sf^m$. 

Recall that the Ricci curvature is nonnegative and at every point there exists a tangent direction 
where the Ricci curvature vanishes. Using that the hypersurface $f$ has a principal curvature of 
multiplicity $n-1$, it follows form the Gauss equation that $M^n$ has nonnegative sectional curvature. 
Given that $M^n$ is not homeomorphic to $\Sf^n$ and oriented, the main result in \cite{fz} yields 
that $M^n$ is diffeomorphic to $\Sf^1\times \Sf^{n-1}$. It follows from part (iii) that 
$f(M^n)$ is isometric to a torus $\mathbb T^n_1(r)$ with $r>1/\sqrt{n}$. 
\qed

\medskip

{\noindent{\it Proof of Corollary \ref{cor1}.}}
Observe that 
\be\label{S}
S\leq 2\sqrt{n-1}\leq a(n,1,H,1).
\ee

At first we claim that if $M^n$ is not a homology sphere, then  $f(M^n)$ is isometric to a torus $\mathbb T^n_1(r)$ for 
appropriate $r$. Suppose that $M^n$ is not a homology sphere. We claim that $b_1(M^n)>0$. Assume to the 
contrary that  $b_1(M^n)=0$, and let $q$ be an integer such that $b_q(M^n)>0$ with $1<q\leq n/2$. Theorem \ref{thm.a}(ii)
implies that $f$ is minimal with $S=n$ and this is a contradiction. Hence $b_1(M^n)>0$ and Theorem \ref{mainthm1}(ii)
shows that equality holds in \eqref{S} pointwise. In particular, we have $H>0$ on $M^n$. The claim now follows
from Theorem \ref{thm.a}(iii).

Now suppose that $f(M^n)$ is not isometric to any torus $\mathbb T^n_1(r)$.
Hence $M^n$ is a homology sphere. 
Using \eqref{ar}, it is obvious that \eqref{S} is written equivalently as
$$
\|\Phi\|^2\leq 2\sqrt{n-1}-nH^2\leq r^2(n,1,H,1), 
$$
where $r(n,1,H,1)$ is the positive root of the polynomial $P$  given by \eqref{P}. Hence we have
$$
P(\|\Phi\|,1,H,1)\geq P\left(2\sqrt{n-1}-nH^2,1,H,1\right)<0
$$
at points where $2\sqrt{n-1}>nH^2$. It follows from inequality \eqref{ineq4} 
for $p=1$ that the Ricci curvature of $M^n$ at these points satisfies
$$
\Ric \geq -\frac{n-1}{n}P\left(2\sqrt{n-1}-nH^2,1,H,1\right)>0.
$$
On the other hand, the points where $2\sqrt{n-1}=nH^2$ are totally umbilical points and consequently the Ricci curvature is 
$$
\Ric = (n-1)(H^2+1).
$$
Ii is clear from the above that the Ricci curvature of $M^n$ is bounded from below by a positive constant. 
The rest of the proof is a consequence of Corollary \ref{cor0}. 
\qed

\medskip

{\noindent{\it Proof of Theorem \ref{thm.A1}.}} 
Since $M^n$ is compact and $f(M^n)$ is contained in an open hemisphere in case where $c=1$, there exists 
a point $x_0\in M^n$ and a unit normal vector $\eta\in N_fM(x_0)$ such that the shape operator 
$A_\eta(x_0)$ is positive definite.

At first we observe that $b_q(M^n)=0$ for every integer $p<q\leq n/2$. Indeed, note that 
$$
S\leq a(n,p,H,c)\leq a(n,q,H,c)
$$
on $M^n$ and $S(x_0)<a(n,q,H(x_0),c)$. Then Theorem \ref{mainthm1}(ii) implies that $b_q(M^n)=0$ 
for every $p<q\leq n/2$.

Now we claim that if $b_p(M^n)>0$, then $p=n/2$ and $M^n$ has the homotopy type of a CW-complex with 
cells only in dimensions $0,n/2$ or $n$. Indeed, if $b_p(M^n)>0$, then 
it follows from Theorem \ref{mainthm1}(iii) that $S= a(n,p,H,c)$
on $M^n$ and that every harmonic form $\omega\in {\mathscr H}^p(M^n)$ is parallel. 
Thus $M^n$ supports a nontrivial parallel $p$-form $\omega$. We distinguish two cases.

{\it Case $c=1$}. Suppose to the contrary that $p<n/2$. We argue on an open neighborhood of $U$ of 
$x_0$ where the mean curvature is positive. Proposition \ref{cod} implies that $f$ is  
a hypersurface in a totally geodesic sphere $\Sf^{n+1}$ of $\Sf^m$ with two distinct principal curvatures 
$\lambda$ and $\mu$ of multiplicities $p$ and $n-p$, respectively, that satisfy $\lambda\mu+1=0$. 
Since the shape operator $A_\eta(x_0)$ is positive definite, obviously $A_\xi(x_0)$ has to be either positive, 
or negative definite. This contradicts the fact that the principal curvatures satisfy $\lambda\mu+1=0$. 
Thus $p=n/2$. That $M^n$ has the homotopy type of a CW-complex with cells only in dimensions $0,n/2$ 
or $n$ follows from Theorem \ref{thm.a}(i).

{\it Case $c=0$}. Suppose to the contrary that $p<n/2$. Equation \eqref{boch.formula} yields 
$\<\B^{[p]}\omega,\omega \>=0$. Hence equality holds in inequality \eqref{ineq4} at every point. 
Let $U$ be an open neighborhood of $x_0$ where the mean curvature is positive. 
Hereafter we argue on $U$. Proposition \ref{lem3}(ii) yields that 
$$
\alpha_f(X,Y)=\<A_\xi X,Y\>\xi \;\;\text{for all} \;\;X,Y\in TU,
$$
where $\xi=\mathcal H/H$ and the shape operator $A_\xi$ has at most two distinct eigenvalues $\lambda$ 
and $\mu$ with multiplicities $p$ and $n-p$, respectively. The Gauss equation is written as
$$
R(X,Y)Z=\<A_\xi Y,Z\>A_\xi X-\<A_\xi X,Z\>A_\xi Y
$$
and consequently 
$$
R(X,Y)\omega=(A_\xi X)^\#\wedge i_{A_\xi Y}\omega-(A_\xi Y)^\#\wedge i_{A_\xi X}\omega.
$$

We choose a local orthonormal frame $X_1,\dots,X_n$ with dual frame $\theta_1,\dots,\theta_n$ such that 
$$
A_\xi X_i=\lambda X_i, \; 1\leq i\leq p
\; \; \text{and}\; \; 
 A_\xi X_j=\mu X_j,\; p+1\leq j\leq n.
$$
Since $\omega$ is parallel, we have $R(X_i,X_j)\omega=0$. From the above it follows that
$$
\lambda\mu\Psi_{ij}=0\;\;\text{for} \;\;1\leq i\leq p, \; p+1\leq j\leq n,
$$
where 
$$
\Psi_{ij}=\theta_i\wedge i_{X_j} \omega-\theta_j\wedge i_{X_i}\omega.
$$
Arguing as in the proof of Proposition \ref{cod}, we have $\Psi_{ij}\neq 0$ for all $1\leq i\leq p<j\leq n$. 
Therefore $\lambda\mu=0$. It is clear that the shape operator 
$A_\xi(x_0)$ has to be either positive, or negative definite, since $A_\eta(x_0)$ is positive definite. 
This contradicts the fact that the principal curvatures 
satisfy $\lambda\mu=0$. Thus $p=n/2$. Since equality holds in inequality \eqref{ineq4} at every point, 
Proposition \ref{lem3} implies that the second fundamental form of $f$ has at most two distinct principal 
curvatures both of multiplicity $p$ in each normal direction. That $M^n$ has the homotopy type of a 
CW-complex with cells only in dimensions $0,n/2$ or $n$ follows by using Morse theory, as in the proof 
of Theorem \ref{thm.a}(i).

Hence, we conclude that either $b_i(M^n)=0$ for all $p\leq i\leq n-p$, or $p=n/2$ and $M^n$ has the 
homotopy type of a CW-complex with cells only in dimensions $0,n/2$ or $n$, and this completes the 
proof of part (i).

Next we assume that $p=1$. From the above we have that $M^n$ is a real homology sphere. 
Our pinching condition implies that the right hand side in inequality \eqref{ineq4} is nonnegative. Thus the 
Ricci curvature of $M^n$ is nonnegative. 
We claim that there is a point where all Ricci curvatures are positive. Arguing by contradiction, we 
suppose that at every point there exists a tangent direction where the Ricci curvature vanishes. 
Hence equality holds in \eqref{ineq4} at every point. Using that 
the mean curvature is positive at the point $x_0\in M^n$, it follows from Proposition \ref{lem3}(ii) that 
$$
\alpha_f(X,Y)=\<A_\xi X,Y\>\xi \;\;\text{for all} \;\;X,Y\in T_{x_0}M,
$$
where $\xi\in N_f(M)(x_0)$ is a unit normal vector and the shape operator $A_\xi$ has an 
eigenvalue with multiplicity $n-1$. Given that the shape operator $A_\eta(x_0)$ is positive definite, 
it is obvious that $A_\xi(x_0)$ is either positive, or negative definite.  In either case, it follows from 
the Gauss equation that the Ricci curvature at $x_0$ is positive for each tangent direction, and 
this is a contradiction.

Since the Ricci curvature is nonnegative and there exists a point where all Ricci curvatures 
are positive, it follows from Aubin \cite{A} that $M^n$ admits a metric with positive Ricci curvature. 
Then the Bonnet-Myers theorem implies that the fundamental group of $M^n$ is finite. In analogy 
with the proof of Corollary \ref{cor}, we can prove that $M^n$ is diffeomorphic to a spherical 
space if $n=3$, and its universal cover is homeomorphic to $\Sf^4$ if $n=4$.

Suppose now that $H>0$ on $M^n,p<(n-1)/2$ and $c=0$. Then we have
$$
S\leq a(n,p,H,0)<a(n,p+1,H,0).
$$
From \cite[Th. 1]{Vlachos07}, it follows that $H_i(M^n;\Z)=0$ for each $p< i< n-p$. Since 
$H_{n-p-1}(M^n;\Z)=0$, the universal coefficient theorem for cohomology (see \cite[Cor. 4, p. 244]{sp}) implies that 
$H^{n-p}(M^n;\Z)$ has no torsion and neither does $H_p(M^n;\Z)$ by Poincar\'e duality. 
Since $b_p(M^n)=0$, we obtain $H_p(M^n;\Z)=0$.
\qed

\medskip

{\noindent{\it Proof of Theorem \ref{thm.d}.}} 
We claim that $b_p(M^n)=0$. Suppose, by contradiction, that $b_p(M^n)>0$. 
Theorem \ref{mainthm1}(iii) implies that $S=a(n,p,H,-1)$ on $M^n$ and that
every harmonic form $\omega\in {\mathscr H}^p(M^n)$ is parallel. Let $\omega$ be 
a non vanishing parallel form. It follows from equation \eqref{boch.formula} that 
$\<\B^{[p]}\omega,\omega \>=0$. Our assumption implies that the right hand side in \eqref{ineq4} 
is nonnegative. Thus, equality holds in inequality \eqref{ineq4} at every point. Proposition \ref{cod} 
implies that $f$ is a hypersurface in a totally geodesic $\Hy^{n+1}\subset \Hy^m$ with principal 
curvatures $\lambda$ and $\mu =1/\lambda$ of multiplicities $p$ and $n-p$, respectively. 

If $p>1$, then the principal curvatures remain constant along
each principal curvature distribution. Hence $f$ is isoparametric. Since $M^n$ is compact, 
$f(M^n)$ is a hypersphere. Consequently, we have $\lambda=\mu =1/\lambda$. The Gauss 
equation implies that $M^n$ is flat, and this is a contradiction since the sectional curvature of 
$M^n$ has to be positive. Now suppose that $p=1$. Since $\<\B^{[1]}\omega,\omega \>=0$, 
at each point there is tangent direction where the Ricci curvature of $M^n$ vanishes. The Gauss 
equation then implies that $\lambda=\mu =1/\lambda$, and this is again a contradiction.

Thus $b_p(M^n)=0$. For each $p<i\leq n/2$ we have
$$
S\leq a(n,p,H,-1)\leq a(n,i,H,-1).
$$
Repeating the above argument, we obtain $b_i(M^n)=0$ for every $p<i\leq n/2$. 
\qed

\section{Integral bound for the Bochner operator}

For the proof of Theorem \ref{thm.integral}, we need to recall some well known facts 
on the total curvature and how Morse theory provides restrictions on the Betti numbers.
Let $f\colon M^n\to \R^{n+k}$ be an isometric immersion of a compact Riemannian manifold into 
the Euclidean space $\R^{n+k}$ equipped with the usual inner product $\langle\cdot,\cdot\rangle$. 
The unit normal bundle is defined by 
$$
UN_f=\left\{(x,\eta)\in N_fM:\|\eta\|=1\right\}.
$$

The {\it generalized Gauss map} $\nu\colon UN_f\to \Sf^{n+k-1}$ is 
given by $\nu(x,\eta)=\eta$. For each $u\in\Sf^{n+k-1}$, we consider the height function 
$h_u$ defined by $h_u(x)=\langle f(x),u\rangle,\, x\in M^n$. 
Since $h_u$ has a degenerate critical point if and only if $u$ is a critical point of $\nu$, by Sard's theorem 
there exists a subset $E\subset\Sf^{n+k-1}$ of measure zero such that $h_u$ is a Morse function for all $u\in\Sf^{n+k-1}\smallsetminus E$. We denote by $\mu_i(u)$ the number of critical 
points of $h_u$ of index $i$ for each $u\in \Sf^{n+k-1}\smallsetminus E$ and set $\mu_{i}(u)=0$ 
for every $u\in E$. Following Kuiper \cite{Kuiper}, we define the {\it total curvature of index $i$ of $f$} by 
$$
\tau_i(f)=\frac{1}{\mathrm{Vol}(\Sf^{n+k-1})}\int_{\Sf^{n+k-1}}\mu_i(u) dS,
$$
where $dS$ denotes the volume element of the sphere $\Sf^{n+k-1}$.
From the weak Morse inequalities (cf. \cite{Milnor63}), we have 
$$
\mu_i(u)\geq b_i(M^n;\Fi) \:\:\text{for all}\;\; u\in \Sf^{n+k-1}\smallsetminus E.
$$ 
By integrating over $\Sf^{n+k-1}$, we obtain 
\begin{equation}\label{TC}
\tau_i(f)\geq b_i(M^n;\Fi).
\end{equation}

There is a natural volume element $d\Sigma$ on the unit normal bundle $UN_f$. In fact, if $dV$ is a 
$(k-1)$-form on $UN_f$ such that its restriction to a fiber of the unit normal bundle at $(x,\eta)$ 
is the volume element of the unit $(k-1)$-sphere of the normal space of $f$ at $x$,
then $d\Sigma=dM\wedge dV$, where $dM$ is the volume element of $M^n$. 
Shiohama and Xu \cite[p. 381]{SX} proved that 
\begin{equation}\label{ShXu}
\int_{U^iN_f}\left|\mathrm{det}A_\eta \right| d\Sigma=\int_{\Sf^{n+k-1}}\mu_i(u) dS,
\end{equation}
where $U^iN_f$ is the subset of the unit normal bundle of $f$ defined by
$$
U^iN_f=\left\{(x,\eta)\in UN_f:\mathrm{Index}\, A_\eta=i\right\},\; \; 0\leq i\leq n.
$$

\medskip

{\noindent{\it Proof of Theorem \ref{thm.integral}.}}
Let $f\colon M^n\rightarrow \R^{n+k}$ be an isometric immersion with second fundamental 
form $\alpha_f$ and shape operator $A_\eta$ with respect to $\eta$, where $(x,\eta)\in UN_f$. 
Using \eqref{bosplit}, \eqref{bo2} and Proposition \ref{int.lemma}, we have
\bea
\big\vert\varrho_p
- F_{n,p}(H, \|\Phi\|)\big\vert^{n/2}
\geq 
(\varepsilon(n,k,p))^{n/2}\int_{\Lambda_p(\alpha_f(x))}\left\vert\det A_\eta \right\vert dV_\eta
\eea
for every point $x\in M^n$. Integrating over $M^n$ and using (\ref{ShXu}), we obtain 
\bea
\int_{M^n} \big| \varrho_p -F_{n,p}(H, \|\Phi\|) \big|^{n/2} dM
\geq\left(\varepsilon(n,k,p)\right)^{n/2} \Vol(\Sf^{n+k-1})\sum_{i=p+1}^{n-p-1} \tau_i(f). 
\eea
Using (\ref{TC}) it follows that 
\begin{eqnarray}\label{57667}
\int_{M^n} \big| \varrho_p -F_{n,p}(H, \|\Phi\|) \big|^{n/2} dM 
\!\!\!&\geq&\!\!\! c(n,k,p)\sum_{i=p+1}^{n-p-1} \tau_i(f) \nonumber \\
\!\!\!&\geq&\!\!\! c(n,k,p)\sum_{i=p+1}^{n-p-1} b_i(M;\Fi), \label{57667}
\end{eqnarray}
where $c(n,k,p)=\left(\varepsilon(n,k,p)\right)^{n/2} \Vol(\Sf^{n+k-1})$.

Now, suppose that
$$
\int_{M^n} \big| \varrho_p -F_{n,p}(H, \|\Phi\|) \big|^{n/2} dM
< c(n,p,k).
$$
It follows from $(\ref{57667})$ that 
$$
\sum_{i=p+1}^{n-p-1} \tau_i(f)<1.
$$ 
Thus there exists 
$u\in \Sf^{n+k-1}$ such that the height function $h_u$ is a Morse function whose number of 
critical points of index $i$ satisfies $\mu_i(u)=0$ for every $p+1\leq i\leq n-p-1$. By Morse theory, 
the manifold $M^n$ has the homotopy type of a CW-complex with no cells of dimension 
$p+1\leq i\leq n-p-1$. 

If $p=1$, there are no $2$-cells and thus by the cellular approximation theorem the inclusion 
of the $1$-skeleton $\mathrm{X}^{(1)}\hookrightarrow M^n$ induces isomorphism between 
the fundamental groups. Therefore, $\pi_1(M^n)$ is a free group on 
$b_1(M^n;\Z)$ elements and $H_1(M^n;\Z)$ is a free abelian group on $b_1(M^n;\Z)$ generators. 
If $\pi_1(M^n)$ is finite, then $\pi_1(M^n)=0$ and hence $H_1(M^n;\Z)=0$. 
From Poincar\'e duality and the universal coefficient theorem, it follows that $H_{n-1}(M^n;\Z)=0$. 
Thus $M^n$ is a simply connected homology sphere and hence a homotopy sphere. 
By the generalized Poincar\'e conjecture (Smale $n\geq 5$, Freedman $n=4$) $M^n$ is 
homeomorphic to $\Sf^n$. \qed

\section{Concluding remarks}
(i) Since homology with real coefficients encodes less topological information that homology with integers ones, it should be interesting to know whether our homology vanishing results hold for  homology with integers coefficients.

(ii) It is worth noticing that the sharp estimate of the Bochner operator in terms of the second fundamental 
form of a submanifold (see Proposition \ref{lem3}) can be used to obtain homology vanishing results for submanifolds in 
different context. For instance, one can strengthen the results in \cite{CalVit}, concerning minimal 
submanifolds in balls with free boundary, by dropping the assumption therein on the flatness 
of the normal bundle.

(iii) DeTurck and Ziller \cite{dtz} obtained minimal isometric embeddings of 
spherical space forms into spheres. None of these minimal submanifolds satisfies the condition $S\leq n$. 
In view of Theorem \ref{thm.b}, this observation raises the question whether any minimal isometric immersion of spherical space forms into spheres is totally geodesic provided that $S\leq n$.

\end{document}